\documentclass[reqno,11pt]{article}

\usepackage{enumerate}
\usepackage{amsmath}
\usepackage{mathrsfs}
\usepackage{amssymb}
\usepackage[round]{natbib}
\usepackage{setspace}
\usepackage{pst-all}
\usepackage{enumerate}
\usepackage[active]{srcltx}
\usepackage{url}
\usepackage{lscape}
\usepackage[bookmarks=true,linkcolor=red,citecolor=blue,urlcolor=green,colorlinks=true, breaklinks]{hyperref}
\usepackage{booktabs}
\usepackage{bigstrut}
\usepackage[refpage]{nomencl} 
\usepackage{longtable}
\usepackage{setspace}
\usepackage{mathtools}
\mathtoolsset{showonlyrefs,showmanualtags}

\newcommand{\rev}[1]{#1}
\newcommand{\revb}[1]{{#1}}

\usepackage{natbib}
 \bibpunct[, ]{(}{)}{,}{a}{}{,}%
 \def\newblock{\ }%

\DeclareMathOperator*{\argmin}{argmin}

\makenomenclature
\newtheorem{thm}{Theorem}
\newtheorem{cor}{Corollary}
\newtheorem{lem}{Lemma}

\newtheorem{dom}{Dominance}

\definecolor{LightOrange} {cmyk}{0,0.7,0.9,0}

\newcommand{\ov}[1]{\overline{#1}}

\newcommand{\lc}[1]{\mathscr{#1}}

\newcommand{\mb}[1]{\mathbf{#1}}
\newcommand{\gb}[1]{\text{\boldmath$#1$}}
\newcommand{\mml}{\{}
\newcommand{\mmr}{\}}
\newcommand{\Real}{\mathbb R}

\newcommand{\lcb}[1]{\overline{\mathscr{#1}}}

\newcommand{\vrpfof}{{\sc{vrpfo}}}
\newcommand{\cvrp}{{\sc{cvrp}}}
\newcommand{\vrp}{{\sc{vrp}}}
\newcommand{\fo}{{\sc{fo}}}
\newcommand{\fp}{{\sc{fp}}}
\newcommand{\lr}{{\sc{lr}}}
\newcommand{\irp}{{\sc{irp}}}
\newcommand{\tcg}{{\sc{cg}}}
\newcommand{\dynp}{{\sc{dp}}}
\newcommand{\fsp}{{\sc{sp}}}
\newcommand{\msp}{F}
\newcommand{\ccg}{{\texttt{DA}}}
\newcommand{\cb}{{\texttt{CB}}}
\newcommand{\deltap}[1]{\Delta(#1)}
\newcommand{\setv}{V}              
\newcommand{\seta}{A}               
\newcommand{\setf}{V_1}             
\newcommand{\setc}{V_2}             
\newcommand{\nf}{n_1}               
\newcommand{\nc}{n_2}               
\newcommand{\nfc}{n}
\newcommand{\setfc}{V_c}            
\newcommand{\st}[1]{s_#1}           
\newcommand{\cij}[2]{d_{#1#2}}      
\newcommand{\tij}[2]{t_{#1#2}}      
\newcommand{\nveh}{m}               
\newcommand{\tveh}{T}               
\newcommand{\rtr}{R}
\newcommand{\cusrtr}{V(R)}
\newcommand{\mmin}{m_{min}}
\newcommand{\mmax}{m_{max}}
\newcommand{\tmin}{\ov{\tveh}}

\newcommand{\dmax}{\Delta^{max}}
\newcommand{\itmax}{\textit{itermax}}
\newcommand{\gmax}{\textit{gapmax}}
\newcommand{\gmin}{\textit{gapmin}}
\newcommand{\tlim}{\textit{tlim}}

\newcommand{\lbcb}{DB_{\cb}}
\newcommand{\ubcb}{PB_{\cb}}
\newcommand{\lbbp}{DB_{\ccg}}
\newcommand{\ubbp}{PB_{\ccg}}

\newcommand{\lbnew}{DB_{new}}

\newcommand{\dk}{{\texttt{DK}}}
\newcommand{\cg}{{\texttt{CG}}}
\newcommand{\genr}{{\texttt{genr}}}

\newcommand{\tabname}{Name}             
\newcommand{\tabtime}{Time}             
\newcommand{\tablbcols}{$|\lcb{R}|$}    
\newcommand{\tablbit}{Iter}             

\newcommand{\tabzbest}{$z^\ast$}        
\newcommand{\tabsolm}{$\#r$}	
\newcommand{\tabnfixed}{$\#o$}	
\newcommand{\tabut}{$\%ut$}
\newcommand{\tabavgc}{$ac$}

\newcommand{\tabopt}{{\tiny (a)}}

\newcommand{\tabnoopt}{}

\newcommand{\tabnstatb}{{\tiny (b)}}

\newcommand{\tabmaxcols}{{\tiny (c)}}

\newcommand{\tabipnotopt}{{\tiny (d)}}

\newcommand{\tablb}{$\%B$}
\newcommand{\taboptip}{IP}

\newcommand{\classCA}{$[\min|c/l]$}
\newcommand{\classPA}{$[\max|p/t]$} 
\newcommand{\ep}[1]{\sigma(#1)}
\newcommand{\nstatb}{NSTATB}

\newcommand{\er}[1]{#1}
\newcommand{\vrptw}{{\sc{vrptw}}}

\allowdisplaybreaks

\RequirePackage{ifthen}
 \renewcommand{\nomgroup}[1]{%
 \ifthenelse{\equal{#1}{A}}{\item[\textbf{A: problem description}]}{%
 \ifthenelse{\equal{#1}{B}}{\item[\textbf{B: mathematical formulation}]}{}
 \ifthenelse{\equal{#1}{C}}{\item[\textbf{C: Lower bounds}]}{}
 \ifthenelse{\equal{#1}{D}}{\item[\textbf{D: Primal bounds}]}{}
 \ifthenelse{\equal{#1}{E}}{\item[\textbf{E: Exact method}]}{}
  }
 }

\begin{document}
\singlespacing

\begin{center}
{\bf\Large Optimal Solution of Vehicle Routing Problems with Fractional Objective Function}
\end{center}
\begin{center}
Roberto Baldacci\\
{\textit{{\small{Department of Electrical, Electronic, and Information Engineering ``Guglielmo Marconi'', University of Bologna, Via Venezia 52, 47521 Cesena, Italy}}}}
\end{center}
\begin{center}
Andrew Lim\\
{\textit{{\small Department of Industrial \& Systems Engineering, National University of Singapore, 1 Engineering Drive 2, Singapore 117576, Singapore}}}
\end{center}
\begin{center}
Emiliano Traversi\\
{\textit{{\small Laboratoire d’Informatique de Paris Nord, Universit\'e de Paris 13; and Sorbonne Paris Cité,
CNRS (UMR 7538), 93430 Villetaneuse, France}}}
\end{center}
\begin{center}
Roberto Wolfler Calvo\\
{\textit{{\small Laboratoire d’Informatique de Paris Nord, Universit\'e de Paris 13; and Sorbonne Paris Cité,
CNRS (UMR 7538), 93430 Villetaneuse, France}}}
\end{center}

\begin{abstract}
 This work proposes a first extensive analysis of the Vehicle Routing Problem with Fractional Objective Function (\vrpfof).
We investigate how the principal techniques used either in the context of fractional programming or in the context of vehicle routing problems interact.
We present new dual and primal bounding procedures which have been incorporated in an exact method.
The method proposed allows to extend specific variants of \vrp\ to their counterpart with a fractional objective function.
Extensive numerical experiments prove the validity of our approach.
\end{abstract}


\onehalfspacing

\section{Introduction}\label{sec:introduction}

In this paper, we investigate how to solve to optimality the Vehicle Routing Problem (\vrp) with optional customers and Fractional Objective Function (\fo). Given a fleet of vehicles and a set of customers, the \vrp\ aims at serving
the customers with a set of feasible routes at minimal cost. \vrp s represent a wide area of combinatorial optimization and mathematical programming. \revb{Most of the works in the literature address the case where the objective
function is linear and all customers must be visited \citep[see][]{TV14}}. In the case considered in this paper, a subset of the customers is \textit{optional}. Therefore, an optional customer is visited only if the objective function improves.
With a linear objective function minimizing the overall cost, the optional customers would never be  part of an optimal solution. However, this is not true when the objective function analyzed is fractional, i.e., it is of the form
$\frac{f(x)}{g(x)}$. In this work we focus on the case where $f(x)$  and $g(x)$ are linear.
Such type of objective function is used to model the so-called \emph{Logistic Ratio} (\lr), which is the ratio of the total cost to the overall resources spent to serve the customers.

In the past, the \lr\ has been widely used in the context of inventory control and production planning (\cite{Bazsa20012} or \cite{Barros1997}). The \lr\ has been also applied to routing problems, more precisely, to the Inventory
Routing Problem (\irp) (\revb{\cite{Morton2011}, \cite{BenoistGJE11}, \cite{GaraixAFJ11} and \cite{Archetti2016}}).

The \irp\ can be described as the combination of a \vrp\ with an inventory management problem, where a supplier has to deliver products to a number of customers over a given time horizon without running out of stock. Therefore,
the \irp\ boils down to  simultaneously decide the inventory management, the vehicle routing and the delivery scheduling (\cite{Coelho2014}).

When an \er{inventory} routing problem is \er{used in practice}, a time horizon needs to be fixed.
Nevertheless, \er{a fixed time horizon, neglects the fact that it would be necessary to solve another optimization problem over the next time horizon.}
This is due to the fact that a direct minimization of costs would lead to postponing as many deliveries as possible to later planning periods (\cite{Dror87}).
It is difficult to define appropriate models able to provide solutions that are ``robust'' also beyond the time horizon considered.
\rev{The \lr\ \er{uses more efficiently} the resources available in \er{the current planning period}, since it \er{provides a solution with better average cost per quantity of
resources used.}}
\er{The needs of customers are better anticipates by minimizing the \lr\, since it avoids myopic behaviours of the optimization models at the end of the time horizon.
While it is not the case with the objective function that directly minimize the cost.}
\revb{What it is more, as observed by \citet{Campbell2001}, real-life inventory routing problems are stochastic. Therefore, any distribution plan covering more than a couple of days
will never be executed completely
as planned. Actual volumes delivered differ from planned volumes because usage rates deviate from their forecasts, \er{the} planned driving time is off due to traffic congestion,
and so forth. Therefore, any planning system needs to be flexible.
In addition, it needs to take advantage of the latest changes in the data, such as last minute orders. Therefore, in practice, long-term plan are often implemented on a daily basis,
and the daily planning needs to capture the costs
and benefits of delivering to a customer earlier than \er{strictly} necessary.}

The \vrpfof\ presents the advantages of planning the current delivery by also taking into account the future demand, while optimizing only over a single day (allowing to keep the
size of the problem under control).
In this context, a customer whose inventory will reach the safety stock level at the end of the day is considered as a compulsory customer (also called \emph{must-go} customer).
On the other hand, a customer whose inventory is below the earliest delivery level but will not reach the safety stock level in the planning window becomes an optional customer
(also called \emph{may-go} customer), see \citet{ortec2016} for more details.
\rev{
Therefore, this work answers to the research question of whether it is possible to find a model able to capture the multiperiod optimization aspect typical of an \irp\ with \fo\ that is also able to handle instances of size of practical interest.}

\subsection{Literature review}\label{sec:literature review}


\rev{The use of the \lr\ as objective function characterizes the problem to solve as a Fractional Programming (\fp) problem, which is a generalization of linear programming.}
\revb{We refer to \cite{Radzik1999}, \cite{Frenk2005}, \cite{schaible2004recent} or \cite{stancu2012fractional} for a brief introduction to \fp.}

Two main techniques are available in the literature to tackle \fp: the variable substitution presented by Charnes and Cooper (\cite{CC62}) and the Dinkelbach's algorithm (\cite{Dinkelbach1967}).
In the former method, the variable substitution allows to rewrite the problem as linear, if the  divisor of the objective function is always greater than zero. The main drawbacks of this technique are that it can potentially lead to numerical instability and that it is difficult to efficiently use it when solving a fractional problem with discrete variables.
Dinkelbach presents an algorithm for problems with (convex) fractional objective function. The algorithm is valid both in case of linear or nonlinear terms in the numerator and denominator. The basic idea of the algorithm is to iteratively solve a parametric linearization of the problem where the parameter represents an estimation of the original objective function. The author also shows that the algorithm terminates in a finite number of iterations.

\revb{A generalization of Dinkelback's algorithm  to nonconvex (continuous) objective functions is presented in~\cite{Rodenas1999}. \citeauthor{Rodenas1999} also extended the Dinkelbach approach to integer linear fractional programming and observed that the algorithm converges on a finite number of iterations.}
\revb{In~\cite{Espinoza2010}, the authors show how lifting, tilting and fractional programming can be viewed as the same optimization problem. \citeauthor{Espinoza2010} described an exact algorithm for the case of mixed integer linear fractional programming and, by combining results of \citet{Dinkelbach1967} and \citet{Schaible1976}, provided a proof of its superlinear convergence rate.
Methods for mixed integer linear fractional programming
with application in cyclic process scheduling problems were also considered by \cite{You2009}. In particular, \citeauthor{You2009} extended Dinkelbach approach to discrete problems by also shoving its superlinear convergence rate.}

\revb{\vrp s with fractional objective functions have been studied by \cite{BenoistGJE11}, \cite{GaraixAFJ11} and \cite{Archetti2016}. \cite{Archetti2016} presented a study on \irp\ with \lr\ and an exact method for its solution. One of the main contributions of their work is a comparison of the optimal solutions obtained when minimizing the \lr\ to the ones obtained when minimizing the total cost.  The solution technique used is an adaptation of the one proposed in~\cite{Dinkelbach1967} to discrete problems. The authors were able to solve to optimality instances involving up to 5 vehicles, 15 customers and 3 periods.
\cite{GaraixAFJ11} investigated the maximization of passenger occupancy rate in \er{a} dial-a-ride problem. The specific objective function considered is to maximize the rate defined as  the sum of the passenger
travel times divided by the total travel time of vehicles. \cite{GaraixAFJ11} proposed two approaches for solving the continuous relaxation of a set partitioning model based on the  Charnes and Cooper's transformation and on Dinkelback's algorithm.
\citet{BenoistGJE11} described a randomized local search algorithm for solving a real-life  routing and  scheduling  problem  arising in optimizing  the  distribution of fluids by  tank  trucks in the  long  run, with the objective to  minimize  the  \lr.}

The book edited by  \citet{TV14} provides a comprehensive overview of exact and heuristic methods for \vrp s.
\revb{In particular, the chapter by \citet{Archetti2014} reviews \vrp s where the set of customers to serve is not given and a profit is associated with each customer that makes such a customer more or less attractive. In this case, the difference between route profit and cost may be maximized, or the profit or the cost optimized with the other measure bounded in a constraint.}

A well know technique used for solving routing problems is column generation (\tcg). Column generation exploits the Dantzig-Wolfe Reformulation of the flow formulation of the original problem and leads to a formulation with an exponential number of variables (we refer the reader to \cite{desaulniers2006column} for a extensive analysis of such a method).

\subsection{Contributions of this paper}\label{sec:contributions}

\revb{In this paper, we describe an exact method to solve the \vrpfof. We formulate the \vrpfof\ using a Set Partitioning (\fsp) like formulation with a fractional objective function. The \vrpfof\ formulation adopted can be used to model alternative objective functions, such as the \lr\ and the maximization of profit over time.

The exact method combines two bounding procedures, derived from the \fsp\ like formulation, with an extension of Dinkelbach's algorithm for fractional programming to integer programs. More precisely, the bounding procedures are used within a route enumeration scheme \citep[see][]{Baldacci2008} to reduce the number of variables of the integer problems solved at each iterations of the Dinkelbach approach.

\er{The extension of the procedure described is possible thanks to the two following novelties presented in this paper:
\begin{itemize}
\item We present a new linear transformation which allows to use a dual ascent heuristic to solve the master problem. It is alternative to the one presented by Charnes and Cooper (\cite{CC62}).
\item We show how the final dual solution of the new linearization can be used to generate a reduced problem containing only the routes whose reduced costs are smaller than a given threshold.
\end{itemize}
}


We perform extensive computational results on instances derived from the \vrp\ literature with different fractional objective functions. The results obtained show that the proposed method is able to solve instances involving up to 79 customers.
}

 The remainder of this paper is organized as follows. In Section \ref{sec:problem description}, we formally describe the problem addressed in this paper and we present the \fsp\ model. In Section \ref{sec:lower bounds}, we present \rev{dual bounds} based on both continuous and integer relaxations of the \fsp\ model. This section also describes the bounding procedures used to compute the \rev{dual bounds}. Section
 \ref{sec:exm} describes the exact method. Dynamic programming algorithms for generating nonelementary routes and feasible
 and elementary routes are described in Section \ref{sec:ng-genr}. We provide the computational studies in Section \ref{sec:computational results} and concluding remarks in Section \ref{sec:conclusions}.

\section{Problem description and mathematical formulation}\label{sec:problem description}

\nomenclature[aa]{$G=(\setv,A$)}{complete digraph}
\nomenclature[ab]{$\setv=\{0\} \cup \setf \cup \setc$, $\nf$,$\nc$, $\nfc$}{sets}
\nomenclature[ac]{$\setfc= \setf \cup \setc$}{customers}
\nomenclature[ad]{$\st{i}$,$\cij{i}{j}$,$\tij{i}{j}$}{times and costs}
\nomenclature[ae]{$\nveh$, $\tveh$}{vehicle fleet}
\nomenclature[af]{$\rtr$, $\cusrtr$, $A(\rtr)$}{route}

The Vehicle Routing Problem with Fractional Objective Function (\vrpfof) considered in this paper can be described as follows.

A complete digraph $G=(\setv,\seta)$ is given, where the vertex set $\setv$ is partitioned as $\setv=\{0\} \cup \setf \cup \setc$. Vertex $0$ represents the depot, vertex set $\setf=\{1,2,\dots,\nf\}$ represents $\nf$
\textit{mandatory} customers, and vertex set $\setc=\{\nf+1,\dots,\nf+\nc\}$ represents $\nc$ \textit{optional} customers. We denote with $\setfc= \setf \cup \setc$, $\nfc=\nf+\nc$ and we assume that $\nf > 0$.
With each vertex $i \in \setv$ is associated a service time $\st{i} > 0$ (we assume $s_0=0$). With each arc $(i,j) \in \seta$ are associated a travel or routing cost $\cij{i}{j}$ and a travel time $\tij{i}{j} \geq 0$.
At the depot it is based a vehicle fleet composed of a set of $\nveh$ identical vehicles. To each vehicle is associated a maximum working time equal to $\tveh$.

A vehicle route $\rtr = (0,i_1,\ldots,i_r,0)$, with $r \geq 1$, is a simple circuit in $G$ passing through the depot, visiting vertices $\cusrtr = \{i_1,\ldots,i_r \}$, $\cusrtr \subseteq \setfc$, and such that the total \textit{working time} computed as the sum of the total service time of the customers visited and the total travel time of the arcs traversed by the route is less than or equal to $\tveh$, i.e., $\sum_{i \in \cusrtr}s_i + \sum_{(i,j) \in \seta(\rtr)}\tij{i}{j} \leq \tveh$, where $\seta(\rtr)$ is the set of arcs traversed by route $R$.
The \textit{cost} of route $\rtr$ is equal to the sum of the travel costs of the arc set traversed by route $\rtr$, i.e., $\sum_{(i,j) \in \seta(\rtr)}\cij{i}{j}$.

We consider the problem of visiting the mandatory customers and of choosing a subset of the optional customers to visit using vehicles based at the depot. More precisely, the \vrpfof\ consists of designing at most $\nveh$ routes such that (i) each vehicle is used at most once, (ii) each mandatory customer is visited once, and (iii) each optional customer is visited at most once.


\er{The \vrpfof\ models the following fractional linear objective functions of practical interest:}

\begin{enumerate}[(i)]
  \item \textit{Minimization of Cost/Load}. Let $q_i$ be the demand associated with each customer $i \in \setfc$ (we assume $q_0=0$). In addition, a fleet of $\nveh$ identical vehicles of capacity $Q$ is stationed at the depot.
      The \textit{load} of a route $\rtr = (0,i_1,\ldots,i_r,0)$ is equal to the total demand of visited customers, i.e. $\sum_{i \in V(\rtr)}q_i$. The objective is to minimize the ratio of the total travel or routing cost divided by the total load of the routes selected in solution (i.e., the \lr). This \vrp\ can be solved as a \vrpfof\ by setting $\tveh=Q$, $s_i=q_i$, $\forall i \in \setfc$, and $t_{ij}=0$, $\forall (i,j) \in A$.
  \item \textit{Maximization of Profit/Time}. Let $p_i$ be a nonnegative profit associated with each customer $i \in \setfc$ (we assume $p_0=0$). The \textit{profit} of a route is equal to the total profit of the visited customers, i.e., $\sum_{i \in V(\rtr)}p_i$. The objective is to maximize the ratio of the total profit divided by \rev{the total working time} of the routes selected in solution. This \vrp\ can be \rev{solved} as a \vrpfof\ by setting $d_{ij}=-p_j$, $\forall (i,j) \in A$.
\end{enumerate}

\revb{It is worthwhile to mention that the special case of the \vrpfof\ where all customers are mandatory (i.e., $n_2=0$) and the objective function is the minimization of Cost/Load, is the Capacitated \vrp\ (\cvrp), which is  in turn a special case of the \vrp\ with Time Windows (\vrptw).

For the state-of-the-art exact algorithms for deterministic \vrp, we refer readers to \citet{Baldacci_2012,TV14,pecin2017new,pecin2017improved}, among others.}

\subsection{Mathematical formulation}\label{sec:mathematical formulation}

\nomenclature[ba]{$\lc{R}$, $\ell$, $\rtr_\ell$, $\setf(\rtr_\ell)$, $\setc(\rtr_\ell)$, $a_{i\ell}$, $c_\ell$, $w_\ell$}{routes}
\nomenclature[ba]{$x_{\ell}$, $F$}{formulation $F$}
In this section, we model the \vrpfof\ as a \fsp\ problem with side constraints and a linear fractional objective function.

Let $\lc{R}$ be the index set of all routes. Given a route $\ell \in \lc{R}$, we denote with $\rtr_\ell$ the sequence $(i_1=0,i_2,\dots,i_r=0)$ of the vertices visited by the route and with $\setf(\rtr_\ell)$ and $\setc(\rtr_\ell)$ the sets $\setf \cap V(\rtr_\ell)$ and $\setc \cap V(\rtr_\ell)$, respectively.
Let $a_{i\ell}$ be a (0-1) binary coefficient  equal to 1 if node $i \in V(\rtr_\ell)$, 0 otherwise. Given a route $\ell$, we denote with $c_\ell$ and $w_\ell$ the routing cost and the working time of route $\ell$, respectively, computed as $c_\ell=\sum_{(i,j)\in \seta(R_\ell)}d_{ij}$  and $w_\ell=\sum_{i \in V(\rtr_\ell)}s_i+\sum_{(i,j)\in \seta(R_\ell)}t_{ij}$.

Let $x_{\ell}$, $\ell \in \lc{R}$, be a (0-1) binary variable equal to $1$ if and only if route $\ell$ is in the optimal solution.
The \vrpfof\ formulation based  on the \fsp\ model, hereafter called $\msp$, is
\begin{alignat}{3}
(\msp)\quad  z(\msp)= \min  & \frac{\sum_{ \ell \in \lc{R}} c_\ell x_\ell}{\sum_{ \ell \in \lc{R}} w_\ell x_\ell } \label{F.objective}\\
                    s.t. & \sum_{\ell \in \lc{R}} a_{i\ell} x_\ell = 1, & \quad \forall i \in \setf             \label{F.degree-f} \\
                         & \sum_{\ell \in \lc{R}} a_{i\ell}x_\ell \leq 1,  & \forall i \in \setc       \label{F.degree-c}\\
                         & \sum_{\ell \in \lc{R}} x_\ell \leq m, \label{F.degree-m}\\
                         & x_\ell \in \{0,1\}, \quad &\forall \ell \in \lc{R}.                       \label{F.int-xi}
\end{alignat}


\er{In the formulation, the objective function states
either
to minimize the \textit{Cost/Load} ratio or
to maximize the \textit{Profit/Time} ratio.}

Constraints \eqref{F.degree-f} and \eqref{F.degree-c} impose that each mandatory customer has to be visited by exactly one route and each optional customer has to be visited at most once by the routes selected in the solution, respectively.
Constraint \eqref{F.degree-m} requires that at most $\nveh$ routes are selected in the solution.


\section{\rev{Dual bounds} for the \vrpfof}\label{sec:lower bounds}

\nomenclature[ca]{$X$, $C\msp$, $z(C\msp)$}{feasible set, continuous relaxation}
\nomenclature[cb]{CC\msp, $z(CC\msp)$, $u$, $y$}{Charnes and Cooper (CC)}
\rev{
In this section, we describe relaxations and bounding procedures for the \vrpfof.
We first present two basic dual bounding techniques, the first one, called \dk, solves directly the continuous relaxation of $\msp$.
The second one, called  \cg, solves a reformulation of the continuous relaxation of $\msp$, proposed for the first time in \citet{CC62}.
In section~\ref{sec:dual_ascent}, we describe an alternative transformation to the one presented in \citet{CC62}, such new transformation is used in an advanced dual ascent bounding procedure called \ccg. Finally, in Section \ref{sec:intlb}, we describe a dual bounding procedure \cb\ based on an integer relaxation of formulation $\msp$}.


Let $X=\{\mb{x} \in \Real_+^{|\lc{R}|}: \eqref{F.degree-f}, \eqref{F.degree-c}, \text{ and } \eqref{F.degree-m}\}$. We define
$z(C\msp)= \min \mml \sum_{ \ell \in \lc{R}} c_\ell x_\ell / \sum_{ \ell \in \lc{R}} w_\ell x_\ell:  \mb{x} \in X\mmr$ as the optimal solution cost of the continuous relaxation of formulation $\msp$, called $C\msp$.
\rev{Since the objective function of formulation $C\msp$ is the quotient of linear functions and $X$ is a convex feasible set, the algorithm proposed by \citet{Dinkelbach1967} can also be used to compute $z(C\msp)$ by means of the solution of a sequence of linear programming problems.
In the following, we identify with \dk\ the bounding procedure corresponding to solve $C\msp$ with the algorithm proposed by \citet{Dinkelbach1967} where each linear programming problem is solved by column generation.}

Under the assumption that $X$ is non-empty and bounded, the following transformation, proposed by \citet{CC62}:
 \er{$$u=\frac{1}{\sum_{ \ell \in \lc{R}} w_\ell x_\ell}$$
 $$y_\ell=u x_\ell,~~ \forall \ell \in \lc{R}$$}
 translates formulation $C\msp$ into the equivalent linear program:
 \begin{alignat}{3}
(CC\msp)\quad  z(CC\msp)= \min  & \sum_{ \ell \in \lc{R}} c_\ell y_\ell \label{CF.objective}\\
                    s.t. & \sum_{\ell \in \lc{R}} a_{i\ell} y_\ell = u, & \quad \forall i \in \setf             \label{CF.degree-f} \\
                         & \sum_{\ell \in \lc{R}} a_{i\ell}y_\ell \leq u,  & \forall i \in \setc       \label{CF.degree-c}\\
                         & \sum_{\ell \in \lc{R}} y_\ell \leq mu, \label{CF.degree-m}\\
                         & \sum_{\ell \in \lc{R}} w_\ell y_\ell =1, \label{CF.yu}\\
                         & u \geq 0, \label{CF.u} \\
                         & y_\ell \geq 0, \quad &\forall \ell \in \lc{R}                       \label{CF.x}.
\end{alignat}
\rev{
\er{$CC\msp$ contains an exponential number of variables. In practice, such problems are solved with a column generation procedure (See~\cite{desaulniers2006column}).}
In the following, we identify with \cg\ the dual bounding procedure that solves (CC\msp) via column generation.
The computational results of Section \ref{sec:computational results} reports a comparison between a \dk\ and \cg\ for computing the dual bound $z(C\msp)$.
}
\subsection{\rev{Dual bounding procedure based on an alternative transformation of $C\msp$ -- \ccg}}\label{sec:dual_ascent}

\nomenclature[cc]{$NC\msp$, $z($NC\msp$)$, $\beta$, $\ov{w}_\ell$, $\ov{a}_\ell$, $\ov{b}_\ell$}{new transformation}
\nomenclature[cd]{$DNC\msp$,$v$}{dual of $NC\msp$}
\nomenclature[ce]{$\lambda_i$, $\phi_i$, $\pi_i$}{theorem dual ascent}
\nomenclature[cf]{$DCC\msp$, $z($DCC\msp$)$, $\mu$, $\omega$}{theorem dual transformation}

The alternative transformation is based on the observation that since $\beta=\sum_{i \in \setf}s_i > 0$, then the term at the denominator of objective function \eqref{F.objective} can be rewritten as $\beta + \sum_{\ell \in \lc{R}}\ov{w}_\ell x_\ell$, where $\ov{w}_\ell=w_\ell - \sum_{i \in \setf(\rtr_\ell)} s_i$. The following lemma holds.

\begin{lem}\label{lem1}
  Let $\mb{y}=\mb{x}/(\beta + \mb{\ov{w}}^T \mb{x})$. Then the objective function of problem $CF$ can be rewritten as $z(C\msp)=\min \mb{c}^T \mb{y}$. In addition, any inequality $\gb{\alpha}^T \mb{x} \leq \alpha_0$, $\gb{\alpha} \in \Real^{|\lc{R}|}$, $\alpha_0 \in \Real$, can be rewritten as $\ov{\gb{\alpha}}^T \mb{y} \leq \alpha_0$ where $\ov{\gb{\alpha}}=(\beta \alpha_1 + \alpha_0 \ov{w}_1, \beta \alpha_2 + \alpha_0 \ov{w}_2, \dots, \beta \alpha_{|\lc{R}|} + \alpha_0 \ov{w}_{|\lc{R}|})$.
\end{lem}
\textit{Proof.}
It is easy to see that $z(C\msp)=\mb{c}^T\mb{x}/(\beta + \mb{\ov{w}}^T \mb{x})=\mb{c}^T\mb{y}$.
Sice $\beta > 0$, we have
\begin{equation}\label{lem1.a}
\begin{split}
  \gb{\alpha}^T \mb{x} \leq \alpha_0 \Rightarrow \gb{\alpha}^T \mb{x} +\frac{\alpha_0}{\beta}\ov{\mb{w}}^T \mb{x} \leq \frac{\alpha_0}{\beta} (\ov{\mb{w}}^T \mb{x} + \beta) \Rightarrow \\
  \frac{\gb{\alpha}^T \mb{x}}{\ov{\mb{w}}^T \mb{x} + \beta} + \frac{\alpha_0}{\beta}\frac{\ov{\mb{w}}^T \mb{x}}{\ov{\mb{w}}^T \mb{x} + \beta} \leq \frac{\alpha_0}{\beta} \Rightarrow \\
  \beta \gb{\alpha}^T \mb{y} + \alpha_0 \ov{\mb{w}}^T \mb{y} \leq \alpha_0 \Rightarrow (\beta \gb{\alpha}^T + \alpha_0 \ov{\mb{w}}^T) \mb{y} \leq \alpha_0.\,\Box
\end{split}
\end{equation}

Formulation $C\msp$ can now be transformed into the following equivalent linear program:
\begin{alignat}{3}
(NC\msp)\quad  z(NC\msp)= \min  & \sum_{ \ell \in \lc{R}} c_\ell y_\ell \label{NCF_objective}\\
                    s.t. & \sum_{\ell \in \lc{R}} \ov{a}_{i\ell} y_\ell = 1, & \quad \forall i \in \setf \label{NCF_degree-f}\\
                         & \sum_{\ell \in \lc{R}} \ov{a}_{i\ell}y_\ell \leq 1,  & \forall i \in \setc \label{NCF_degree-c}\\
                         & \sum_{\ell \in \lc{R}} \ov{b}_\ell y_\ell \leq m, \label{NCF_degree-m}\\
                         & y_\ell \geq 0, \quad &\forall \ell \in \lc{R}, \label{NCF_x}
\end{alignat}
where $\ov{a}_{i\ell}=\beta a_{i \ell}+\ov{w}_\ell$ and $\ov{b}_\ell = \beta + m \ov{w}_\ell$, $\ell \in \lc{R}$.
Notice that formulation $NC\msp$ has the same number of variables and constraints of formulation $C\msp$. Clearly, $z(NC\msp)=z(CC\msp)=z(C\msp)$.


We denote with $DNC\msp$ the dual of $NC\msp$. The variables of $DNC\msp$ are given by the vector $\mb{v}=(v_0,v_1,\dots,v_{\nfc})$, where $v_1,\dots,v_{\nf} \in \Real$
are associated with constraints \eqref{NCF_degree-f}, $v_{\nf+1},\dots,v_{\nfc} \leq 0$, with constraints \eqref{NCF_degree-c}, and $v_0 \leq 0$ with constraint \eqref{NCF_degree-m}.

\er{Let $DCC\msp$ be the dual of problem $CC\msp$. The variables of $DCC\msp$ are given by vector $\gb{\mu}=(\mu_0,\dots,\mu_{\nfc})$ and variable $\omega$, where $\mu_1, \dots, \mu_{\nf} \in \Real$
are associated with constraints \eqref{CF.degree-f}, $\mu_{\nf+1}, \dots, \mu_{\nfc}  \leq 0$ with constraints \eqref{CF.degree-c}, $\mu_0 \leq 0$ with constraint \eqref{CF.degree-m}, and $\omega \in \Real$
with constraint \eqref{CF.yu}. The following theorem shows how to compute a solution of $DCC\msp$ given a solution of $DNC\msp$.
\begin{thm}\label{thm2}
Let $\mb{v}$ be a feasible solution of problem $DNC\msp$ of cost $z(DNC\msp)$. A feasible solution $(\gb{\mu}, \omega)$ of $DCC\msp$ of cost $z(DCC\msp)=z(DNC\msp)$ can be obtained by setting:
\begin{equation}\label{thm2-a}
  \omega=\sum_{i \in \setfc}v_i + \nveh v_0, \;\; \mu_0=\beta v_0, \;\; \mu_i=\beta v_i - s_i \omega, \forall i \in \setf, \;\; \mu_i=\beta v_i, \forall i \in \setc.
\end{equation}
\end{thm}
\textit{Proof.} The proof is provided in the e-companion to this paper. $\Box$
}

\subsubsection{Dual Ascent Procedure \ccg}
\er{The structure of the new formulation $NC\msp$ allows to use}
a bounding procedure, called \ccg, that is used to compute a near-optimal dual solution of problem $NC\msp$. Procedure \ccg\ differs from standard column generation methods based on the simplex algorithm as it uses a dual
ascent heuristic to solve the master problem \citep[see][]{Baldacci2010a}.

The bounding procedure is based on the following theorem.

\begin{thm}\label{thm1}
Let us associate penalties $\lambda_i \in \Real$, $\forall i \in \setf$, with constraints \eqref{NCF_degree-f}, $\lambda_i \leq 0$, $\forall i \in \setc$, with constraints \eqref{NCF_degree-c}, and $\lambda_0 \leq 0$, with constraint \eqref{NCF_degree-m}.
Let $\lc{R}_i = \{\ell \in \lc{R}: \ov{a}_{i\ell} > 0\}$. For each $i \in \setf$ compute:
\begin{equation}\label{thm1.a}
  \phi_i = \pi_i \min_{\ell \in \lc{R}_i} \left \mml \frac{c_\ell - \lambda(R_\ell) - \ov{b}_\ell \lambda_0}{\pi(R_\ell)} \right \mmr
\end{equation}
where $\pi_i > 0$ is a weight assigned to customer $i \in \setf$, $\lambda(R_\ell)=\sum_{i \in \setfc} \ov{a}_{i \ell}\lambda_i$ and $\pi(R_\ell)=\sum_{i \in \setf}\ov{a}_{i\ell}\pi_i$.
A feasible $DNC\msp$ solution $\mb{v}$  of cost $z(DNC\msp(\gb{\lambda}))$ is given by the following expressions:
\begin{equation}\label{eq:thm1.b}
v_i= \phi_i + \lambda_i, \forall i \in \setf, \quad v_i= \lambda_i, \forall i \in \setc, \quad v_0=\lambda_0.
\end{equation}
\end{thm}
\textit{Proof.} See \citet{Baldacci2010a}. $\Box$

The optimal solution cost of the following problem
\begin{equation}\label{eq:ldual}
 \max_{\mb{\lambda}}\mml z(DNC\msp(\gb{\lambda}))\mmr
\end{equation}
provides the best possible \rev{dual bound} which can be computed
by means of Theorem \ref{thm1}. In practice, problem \eqref{eq:ldual}
cannot be solved even by means of subgradient optimization
as the computation of solution
$\mb{v}$, for given vector $\gb{\lambda}$ requires the a priori generation of the set $\lc{R}$.
Method \ccg\ is an iterative algorithm which computes
a \rev{dual bound} as the cost of a suboptimal solution of problem \eqref{eq:ldual} by using a limited subset $\lcb{R} \subseteq \lc{R}$ and by changing the values of
vector $\gb{\lambda}$.
At each iteration, \ccg\ uses expressions
\eqref{eq:thm1.b} to find a solution $\mb{v}$, for given $\gb{\lambda}$, of the reduced $DNC\msp$ problem defined on route subset $\lcb{R}$ instead of $\lc{R}$. A pricing
procedure is used to
identify the route subset $\lc{N} \subset \lc{R} \setminus \lcb{R}$ whose dual constraints
are violated by the current solution $\mb{v}$. In case $\lc{N} \neq \emptyset$, then $\mb{v}$ is not a feasible $DNC\msp$ solution, and $\lc{N}$ is added to the current core
problem $\lcb{R}$.
At each iteration, subgradient vectors are computed and
used to change vector $\gb{\lambda}$ to maximize the value of the \rev{dual bound}.

Section \ref{sec:ngroute} describes the method used to compute set $\lc{N}$ - for the details of method \ccg\, the reader is referred to \citet{Baldacci2010a}.

\subsection{\rev{Dual bounding procedure based on an integer relaxation of $F$ -- \cb}}\label{sec:intlb}

\nomenclature[ci]{$\tmin$, $\lambda_i$, $\lc{R}^t_i$ }{integer relaxation I}
\nomenclature[cl]{$R\msp/\ov{R\msp}(\ov{t},\ov{m},\gb{\lambda})$, $Maxit3$}{integer relaxation II}
\nomenclature[cm]{$\theta_i$, $\epsilon$, $Maxit3$}{integer relaxation III}
\nomenclature[cn]{$\mmin$, $\mmax$}{dual and primal bounds on $\nveh$}

In this section, we describe an \er{alternative} \rev{dual bound} derived from an integer relaxation of formulation $\msp$. The relaxation is based on the observation that the optimal solution cost $z(\msp)$ of formulation $\msp$ can be computed as
  $z(\msp) = \min_{{\substack{ \tmin \leq \ov{t} \leq \nveh \tveh \\ \left \lceil \ov{t}/\tveh \right \rceil \leq \ov{m} \leq \nveh}}} \left \mml z(\msp(\ov{t},\ov{m}))/\ov{t} \right \mmr$
, where
\begin{alignat}{3}
(\msp(\ov{t},\ov{m}))\quad  z(\msp(\ov{t},\ov{m}))= \min  & \sum_{ \ell \in \lc{R}} c_\ell x_\ell \label{FP.objective}\\
                    s.t. & \eqref{F.degree-f}, \eqref{F.degree-c} \text{ and }\\
                         & \sum_{\ell \in \lc{R}} x_\ell = \ov{m}, \label{FP.degree-m}\\
                         & \sum_{\ell \in \lc{R}}w_\ell x_\ell = \ov{t} \label{FP.t}\\
                         & x_\ell \in \{0,1\}, \quad &\forall \ell \in \lc{R},                       \label{FP.int-x}
\end{alignat}
and $\tmin$ is a valid \rev{dual bound} on the total working time of any feasible solution of formulation $\msp$.
In practice, problem $\msp(\ov{t},\ov{m})$ cannot be solved directly but a valid \rev{dual bound} on its optimal solution cost $z(\msp(\ov{t},\ov{m}))$ can be obtained as follows.

Let $\lambda_i$, $\forall i \in \setfc$, be a set of Lagrangian penalties associated with constraints \eqref{F.degree-f} and \eqref{F.degree-c}, where
$\lambda_i \in \Real$, $\forall i \in \setf$, and $\lambda_i \leq 0$, $\forall i \in \setc$. The Lagrangian relaxation of formulation $\msp(\ov{t},\ov{m})$ by means of penalty vector $\gb{\lambda}$, is as follows:
\begin{alignat}{3}
(R\msp(\ov{t},\ov{m},\gb{\lambda}))\quad  z(\msp(\ov{t},\ov{m},\gb{\lambda}))= & \min \sum_{\ell \in \lc{R}}\ov{c}_\ell x_\ell + \sum_{i \in \setfc}\lambda_i \label{RFP.objective}\\
    s.t. & \eqref{FP.degree-m}, \eqref{FP.t} \text{ and }\eqref{FP.int-x},
\end{alignat}
where $\ov{c}_\ell=c_\ell - \lambda(R_\ell)$  with $\lambda(R_\ell)=\sum_{i \in \setfc}a_{i\ell}\lambda_i$.
Let $\lc{R}^t_i \subseteq \lc{R}$  be the index set of all routes in $\lc{R}$ ending in $i \in \setfc$ and with a total working time $t$. We have $\lc{R}=\bigcup_{\substack{i \in \setfc \\ 1 \leq t \leq T}} \lc{R}^t_i$ and $\lc{R}^t_i \cap \lc{R}^{t'}_j = \emptyset$, $\forall i,j \in \setfc$, $i \neq j$, $t,t' \in [1,\tveh]$.

Let $\varphi^t_i$, $i \in \setfc$, $1 \leq t \leq \tveh$, be a \rev{dual bound} on the modified cost $\ov{c}_\ell$ of any route $\ell \in \lc{R}^t_i$, i.e, $\varphi^t_i \leq \ov{c}_\ell$, $\forall \ell \in \lc{R}^t_i$. We have
\begin{equation}\label{eq:cb-lagr}
  \sum_{\ell \in \lc{R}}\ov{c}_\ell x_\ell = \sum_{i \in \setfc}\sum_{t=1}^{\tveh}\sum_{\ell \in \lc{R}^t_i}\ov{c}_\ell x_\ell \geq \sum_{i \in \setfc} \sum_{t=1}^{\tveh} \varphi^t_i \sum_{\ell \in \lc{R}^t_i}x_\ell.
\end{equation}
We assume $\varphi^t_i=\infty$, $i \in \setfc$, if no feasible routes ending in $i \in \setfc$ with a total working time equal to $t$ exists.
By setting $z^t_i=\sum_{\ell \in \lc{R}^t_i}x_\ell$, from relaxation $R\msp(\ov{t},\ov{m},\gb{\lambda})$ we obtain the following relaxation:
\begin{alignat}{3}
(\ov{R\msp}(\ov{t},\ov{m},\gb{\lambda}))\quad  z(\ov{R\msp}(\ov{t},\ov{m},\gb{\lambda}))= & \min \sum_{i \in \setfc}\sum_{t=1}^{\tveh} \varphi^t_i z^t_i \label{RFP.objective}\\
s.t. & \sum_{i \in \setfc}\sum_{t=1}^{\tveh}  t z^t_i = \ov{t},
 \label{RFP.time} \\
     & \sum_{i \in \setfc}\sum_{t=1}^{\tveh} z^t_i = \ov{m}, \label{RFP.vehicles} \\
     & \sum_{t=1}^{\tveh} z^t_i \leq 1, & \forall i \in \setfc \label{RFP.setfc} \\
     & z^t_i \in \{0,1\}, & \quad \forall i \in \setfc, \forall t, 1 \leq t \leq \tveh.
\end{alignat}
Problem $R\msp(\ov{t},\ov{m},\gb{\lambda})$ can be solved by \dynp\ as follows.
Let $g_i(t,k)$ be the optimal solution to problem $R\msp(\ov{t},\ov{m},\gb{\lambda})$ with the right-hand-side of equation \eqref{RFP.time} replaced by $t$, the right-hand-side of equation \eqref{RFP.vehicles} replaced by $k$, and with $z^t_j=0$ for $j > i$, $t=1,\dots,\tveh$. Function $g_i(t,k)$ can be computed as follows:
\begin{equation}\label{eq:cb-dyng}
  g_i(t,k) = \min \left \mml g_{i-1}(t,k), \min_{1 \leq t' \leq \min \mml t-1, \tveh \mmr} \mml g_{i-1}(t-t',k-1) + \varphi^{t'}_i \mmr \right \mmr
\end{equation}
for $k=2,\dots,\ov{m}$, $i=k,\dots,\nfc-\ov{m}+k$, and $\forall t$, $\max \mml 1, \tmin-(\ov{m}-k)T \mmr \leq t \leq \min \mml k \tveh, \ov{t} \mmr $.
The following initialization is required:
  $g_i(t,1) = \min \mml g_{i-1}(t,1), \varphi^t_i \mmr, \; i=2,\dots,\nfc, \; \max \mml 1, \tmin-(\ov{m}-1)T \mmr \leq t \leq \min \mml \tveh, \ov{t} \mmr$
and
$g_1(t,1) = \varphi^t_i$, $\max \mml 1, \tmin-(\ov{m}-1)T \mmr \leq t \leq \min \mml \tveh, \ov{t} \mmr$, \; $g_i(t,i+1)=\infty$, $i=1,\dots,\ov{m}-1$, $\forall t \in [0,\ov{t}]$.
The optimal value $z(\ov{R\msp}(\ov{t},\ov{m},\gb{\lambda}))$ can then be computed as
  $z(\ov{R\msp}(\ov{t},\ov{m},\gb{\lambda})) = g_{\nfc}(\ov{t}, \ov{m})$.

\subsubsection{Bounding procedure \cb}\label{sec:procbp}

Based on relaxation $R\msp(\ov{t},\ov{m},\gb{\lambda})$ we designed a bounding procedure, called \cb.
Procedure \cb\ is based on the observation that all functions $g_{\nfc}(\ov{t},\ov{m})$ can be computed using the \dynp\ recursion \eqref{eq:cb-dyng} with $\ov{t}=m\tveh$ and $\ov{m}=m$ and
for $k=2,\dots,\nveh$, $i=k,\dots,\nfc$, and $\forall t$, $\max \mml 1, \tmin-(\nveh-k)T \mmr \leq t \leq \min \mml k \tveh, m\tveh \mmr $.
Further, Procedure \cb\ uses subgradient optimization to maximize the value of \rev{dual bound} $z(R\msp(\ov{t},\ov{m},\gb{\lambda}))$.
\rev{Bounding procedure \cb\ works as follows:}

\begin{enumerate}[Step 1.]
    \item \textit{Initialization}. Initialize the penalty vector $\gb{\lambda}=\mb{0}$. Set $\lbcb=-\infty$, $i=1$, $DB(\ov{t},\ov{m})=-\infty$, $\forall \ov{t}, 0 \leq \ov{t} \leq T$, $\forall \ov{m}, 0 \leq \ov{m} \leq m$.
    \item \textit{Compute functions $\varphi^t_i$}. Compute functions $\varphi^t_i$, $\forall i \in \setfc$, $\forall t, \max \mml 1,\tmin-(\nveh-1)\tveh \mmr \leq t \leq \tveh$ (see Section \ref{sec:ngroute}) and set $\varphi^t_i=\infty$ $\forall i \in \setfc$, $\forall t, t < \max \mml 1,\tmin-(\nveh-1)\tveh \mmr$.
    \item \textit{\rev{Dual bound} computation}. Compute function $g_i(t,k)$ using \dynp\ recursion \eqref{eq:cb-dyng} with $\ov{t}=m \tveh$ and $\ov{m}=\nveh$.
        Compute $DB(\ov{t},\ov{m}) = \max \left \mml DB(\ov{t},\ov{m}), \frac{g_{\nfc}(\ov{t},\ov{m})}{\ov{t}} \right \mmr $, $\forall \ov{t}, 1 \leq \ov{t} \leq T$, $\forall \ov{m}, 1 \leq \ov{m} \leq m$.
        Compute
            \begin{equation}\label{eq:cb-dynp-prob}
              z^\ast = \min_{{\substack{ \tmin \leq \ov{t} \leq \nveh \tveh \\ \left \lceil \ov{t}/\tveh \right \rceil \leq \ov{m} \leq \nveh}}} \left \mml DB(\ov{t},\ov{m}) \right \mmr
            \end{equation}
            and let $t^\ast$ and $m^\ast$ be the values producing $z^\ast$ in expression \eqref{eq:cb-dynp-prob}.
            If $z^\ast > \lbcb$, set $\lbcb=z^\ast$.
    \item \textit{Update the penalty vector $\gb{\lambda}$}. Compute the solution corresponding to \rev{dual bound} $z^\ast$ by backtracking
        using recursions \eqref{eq:cb-dyng} and \eqref{ng.a11} and values $t^\ast$ and $m^\ast$. Let $\lcb{R}$ be the set of $(NG,t,i)$-route selected in solution.

        Let $\theta_i$ be the number of times that customer $i \in \setfc$ is visited by the routes in $\lcb{R}$, i.e. $\theta_i=\sum_{\ell \in \lcb{R}}a_{il}$. The value of $\gb{\lambda}$ is modified as follows: $\lambda_i=\lambda_i - \epsilon \gamma (\theta_i-1)$, $\forall i \in \setf$, $\lambda_i=\min \mml 0, \lambda_i - \epsilon \gamma (\theta_i-1) \mmr$, $\forall i \in \setc$, where $\epsilon$ is a positive constant and $\gamma= |0.2 z^\ast| / (\sum_{i \in \setfc}(\theta_i-1)^2)$.
    \item \textit{Termination criteria}.  Set $i = i + 1$. If $i=Maxit3$, stop, otherwise go to Step 2.
\end{enumerate}
At the \rev{end} of the procedure, the value of $\lbcb$ represent the best \rev{dual bound} computed by the procedure. Let $\mmin$ and $\mmax$ be the minimum and maximum values of $\ov{m}$, $\lceil \tmin/\tveh \rceil \leq \ov{m} \leq \nveh$, such that $\min_{\tmin \leq \ov{t} \leq \nveh \tveh} \mml DB(\ov{t}, \ov{m})\mmr < UB$, where $UB$ is a given \rev{primal bound} on the optimal \vrpfof\ solution cost. Values $\mmin$ and $\mmax$ represent valid \rev{dual and primal bound}s on the number of vehicles used in any optimal \vrpfof\ solution.

%
%

\section{An exact method for solving the \vrpfof}\label{sec:exm}

\revb{The \vrpfof\ can be solved to optimality using the extension of the method proposed by \citet{Dinkelbach1967} to integer problem once a generic exact method for solving problem  \eqref{F.degree-c}-\eqref{F.int-xi} with a linear objective function is available; to our knowledge, this problem has never been addressed in the literature. A main drawback of this procedure is that the generic exact method must be applied from scratch at each iteration of the Dinkelbach's algorithm. In this section, we describe an exact method that combines the bounding procedures \ccg\ and \cb\, described in Section \ref{sec:dual_ascent} and Section \ref{sec:intlb}, respectively, to a priori generate a reduced $\msp$ problem containing all routes of any optimal solution, thus reducing the dimensions of the integer problems solved at each iteration of the Dinkelbach's algorithm. The exact method relies on a technique, called \textit{route enumeration}, used by \citet{Baldacci2008} to solve the \cvrp\ and by \citet{Baldacci2011} to solve both the \cvrp\ and the \vrptw. In particular, in Section \ref{sec:varred}, we revisit the technique for fractional linear objectives.}

More precisely, in the exact method we use bounding procedure \ccg\ to
obtain a near-optimal solution $(\gb{\mu}, \omega)$ of the dual problem $DCC\msp$ of cost $z(DCC\msp)$. This solution allows us to compute the reduced costs
  $\ov{c}_\ell = c_\ell - \sum_{i \in \setfc}a_{i\ell}\mu_i - \mu_0 - w_\ell \omega$ of each route $\ell \in \lc{R}$.
Whenever the reduced cost $c_\ell$ of a route $\ell \in \lc{R}$ exceeds a given threshold, computed as a function of a known \rev{primal bound} $UB$ and the dual bound $z(DCC\msp)$, it is possible to eliminate route $\ell$
from $\lc{R}$. Nevertheless, the resulting $\msp$ might still be too large to be solved exactly. We propose an iterative procedure for solving the \vrpfof\ where at each iteration a reduced $\msp$ problem is solved.
Further, each reduced  $\msp$ is solved  to optimality using an exact method based on the extension of the method proposed by \citet{Dinkelbach1967} to integer problem.
The procedure terminates when either an optimal $\msp$ solution is achieved or the distance of the solution cost of the reduced $\msp$ problem from the \rev{dual bound} is less than a user-defined value or a maximum number
of iterations is reached.

\rev{The bounding procedures \ccg\ and \cb\, described in Section \ref{sec:dual_ascent} and Section \ref{sec:intlb} respectively,} are interwoven with a \textit{Lagrangean heuristic} that produces a feasible \vrpfof\ solution.
More specifically, whenever an improved \rev{dual bound} has been computed (see Step 3 of both procedure \ccg\ and \cb), the procedure calls an algorithm that produces a feasible \vrpfof\ solution using the route set $\lcb{R}$.

In the following, we describe the details of the exact method.

%
%

\subsection{Variable reduction of formulation $\msp$}\label{sec:varred}

\rev{The aim of this section is to identify a criterion to restrict ourself to only the column that can be potentially part of the optimal solution, we call such procedure \emph{variable reduction.}}
 Let $(\gb{\mu}, \omega)$ be a feasible solution of $DCC\msp$ of cost $z($DCC\msp$)$ \er{equal to $\omega$, since $\omega$ is the only variable in the objective function with coefficient one.} 
 Let $\ov{\mb{x}}$ be a feasible solution of $\msp$ of cost $\ov{z}(\msp)$ and let $\ov{c}_\ell$ be the reduced cost of route $\ell \in \lc{R}$ with respect to the dual solution $(\gb{\mu}, \omega)$,
 that is $\ov{c}_\ell = c_\ell - \sum_{i \in \setfc}a_{i\ell}\mu_i - \mu_0 - w_\ell \omega$, and let $\ov{c}_0=\sum_{i \in \setfc}\mu_i + \nveh \mu_0$.

The variable reduction of formulation $\msp$ is based on the following Theorem \ref{thm3}:

\begin{thm}\label{thm3} Let $\lcb{R}=\{\ell: \ov{x}_\ell=1, \ell \in \lc{R}\}$. The following inequality holds:
 \begin{equation}\label{eq:thm3-0}
   \ov{z}(\msp) \geq z(DCC\msp) + \frac{\sum_{\ell \in \lcb{R}}\ov{c}_\ell + \ov{c}_0}{\sum_{\ell \in \lcb{R}}w_\ell}.
 \end{equation}
 \end{thm}
\textit{Proof.}
From the definition of $\ov{z}(\msp)$ we have:
\begin{equation}\label{eq:thm3-1}
  \begin{split}
    \ov{z}(\msp)=\frac{\sum_{\ell \in \lcb{R}}c_\ell}{\sum_{\ell \in \lcb{R}}w_\ell} = \frac{\sum_{\ell \in \lcb{R}}(\sum_{i \in \setfc}a_{i\ell}\mu_i+\mu_0+w_\ell\omega + \ov{c}_\ell)-\sum_{i \in \setfc}\mu_i-\nveh\mu_0+\ov{c}_0}{\sum_{\ell \in \lcb{R}}w_\ell}=\\
    \omega \frac{\sum_{\ell \in \lcb{R}}w_\ell }{\sum_{\ell \in \lcb{R}}w_\ell } + \frac{\sum_{i \in \setfc}(\sum_{\ell \in \lcb{R}}a_{i\ell}-1)\mu_i + \sum_{\ell \in \lcb{R}}(1-m)\mu_0 + \sum_{\ell \in \lcb{R}}\ov{c}_\ell + \ov{c}_0}{\sum_{\ell \in \lcb{R}}w_\ell }.
  \end{split}
\end{equation}
Since $\ov{x}$ represents a feasible \vrpfof\ solution, we have (i)  $\sum_{\ell \in \lcb{R}}(a_{i\ell}-1)=0$, $\forall i \in \setf$, (ii)
$\sum_{\ell \in \lcb{R}}(a_{i\ell}-1) \leq 0$, $\forall i \in \setc$, and (iii) $\sum_{\ell \in \lcb{R}}(1-m) \leq 0$, therefore as $\mu_i \leq 0$, $\forall i \in \setc$, and $\mu_0 \leq 0$
\begin{equation}\label{eq:thm3-2}
    \begin{split}
  \sum_{i \in \setfc}(\sum_{\ell \in \lcb{R}}a_{i\ell}-1)\mu_i + \sum_{\ell \in \lcb{R}}(1-m)\mu_0 = \\ \sum_{i \in \setf}(\sum_{\ell \in \lcb{R}}a_{i\ell}-1)\mu_i + \sum_{i \in \setc}(\sum_{\ell \in \lcb{R}}a_{i\ell}-1)\mu_i+\sum_{\ell \in \lcb{R}}(1-m)\mu_0 \geq 0.
    \end{split}
\end{equation}

From equation \eqref{eq:thm3-1} and inequality \eqref{eq:thm3-2} we obtain \eqref{eq:thm3-0}.
$\,\Box$

\begin{cor}\label{cor1}
  Let $UB$ be the cost of a feasible \vrpfof\ solution and let $z($DCC\msp$)$ be the cost of a feasible dual solution $(\gb{\mu}, \omega)$ of $DCC\msp$.
  Any optimal solution $x$ of cost $z(F)$ less than $UB$ cannot contain any route $\ell \in \lc{R}$ such that
      $\ov{c}_\ell \geq \alpha_\ell UB - (\alpha_\ell z(DCC\msp)  +\ov{c}_0)$,
    where $\alpha_\ell=w_\ell + (\nveh-1)\tveh$ \er{is an upper bound on the total working time of any solution containing route $\ell$, since $\nveh$ is the maximum number of vehicles in solution.}
\end{cor}

\textit{Proof.}(By contradiction)
Let $\lcb{R}$ be the index set of the routes of the feasible solution $x$ of cost $z(F) < UB$ and suppose that exists $\ell' \in \lcb{R}$ such that $\ov{c}_{\ell'} \geq \alpha_{\ell'} UB - (\alpha_{\ell'} z(DCC\msp)  +\ov{c}_0)$.

From Theorem \ref{thm3} and as $\ov{c}_\ell \geq 0$, $\forall \ell \in \lcb{R}$, and since $\alpha_{\ell'} \geq \sum_{\ell \in \lcb{R}}w_\ell$ we have:
\begin{equation}\label{eq:cor1-2}
  \begin{split}
   z(\msp) \geq z(DCC\msp) + \frac{\sum_{\ell \in \lcb{R}}\ov{c}_\ell + \ov{c}_0}{\sum_{\ell \in \lcb{R}}w_\ell} \geq z(DCC\msp) + \frac{\ov{c}_{\ell'} + \ov{c}_0}{\sum_{\ell \in \lcb{R}}w_\ell} \geq z(DCC\msp) + \frac{\ov{c}_{\ell'} + \ov{c}_0}{\alpha_{\ell'}} \geq \\ z(DCC\msp) + \frac{\alpha_{\ell'} UB - (\alpha_{\ell'} z(DCC\msp)  +\ov{c}_0) + \ov{c}_0}{\alpha_{\ell'}} \geq UB.\,\Box
  \end{split}
\end{equation}

\rev{
Corollary~\ref{cor1} will be used in the exact procedure  presented in Section~\ref{sec:exmethod} to generate only the columns with a reduced cost lower than the gap $\alpha_\ell UB - (\alpha_\ell z(DCC\msp)  +\ov{c}_0)$.
}

\subsection{Description of the exact method}\label{sec:exmethod}

\nomenclature[ea]{$\dmax$, \gmax, \itmax}{exact method I}
\nomenclature[eb]{$\lbcb$, $\ubcb$, $\lbbp$, $\ubbp$}{exact method II}
\nomenclature[eb]{$\msp(\lcb{R})$, $\gamma_\ell$, $\gmin$,$\lbnew$, $\varepsilon_1$,$\varepsilon_2$}{exact method III}

The method is based on an iterative procedure where at each iteration a reduced version of $\msp$ involving at most $\dmax$ routes ($\dmax$ is a user-defined parameter) is solved.
This procedure terminates when one of the following three conditions is encountered: (i) an optimal $\msp$ solution is achieved, (ii) the distance of the solution cost of the reduced $\msp$ problem from the \rev{dual bound} is less than the user defined value, \gmax, and (iii) the maximum number of iterations, \itmax, is reached.
The scheme of the proposed exact method for solving $\msp$ is as follows.
\begin{enumerate}[Step 1.]
    \item \textit{Initialization}. Set $i=0$ and $z^\ast=\infty$. Initialize \itmax, $\dmax$ and \tlim.
    \item \textit{Execute bounding procedure \cb}. Let $\lbcb$ and $\ubcb$ be the final \rev{dual and primal bounds} computed, respectively. Let $\mmin$ and $\mmax$ be computed as described in Section \ref{sec:procbp}.
    \item \textit{Execute bounding procedure \ccg}. Let $\lbbp$ and $\ubbp$ be the final \rev{dual and primal bounds}  computed and let $(\gb{\mu}, \omega)$ be the solution of problem $DCC\msp$ computed by means of Theorem \ref{thm2} using solution $\ov{\mb{v}}$ corresponding to $\lbbp$.
        Let $\ov{c}_\ell$ be the reduced cost of route $\ell \in \lc{R}$ with respect to the dual solution $(\gb{\mu}, \omega)$, that is $\ov{c}_\ell = c_\ell - \sum_{i \in \setfc}a_{i\ell}\mu_i - \mu_0 - w_\ell \omega$,
        and let $\ov{c}_0=\sum_{i \in \setfc}\mu_i + \nveh \mu_0$.
        Set $z^\ast=\min \mml \ubcb, \ubbp\mmr$ and let $\ov{\mb{x}}$ be the corresponding solution.
    \item \textit{Define a reduced problem $\msp(\lcb{R})$ from $\msp$}. Set $i=i+1$. Generate the largest route set $\lcb{R} \subseteq \lc{R}$ such that:
        \begin{equation}
            \left .
                \begin{array}{ll}
                a)  \; &  |\lcb{R}| \leq \dmax,  \\
                b)  \; &  \ov{c}_\ell < \gamma_\ell, \forall \ell \in \lcb{R},
                \end{array} \hspace{1cm}
            \right \} \label{DynP.Conds}
        \end{equation}
        where $\gamma_\ell=\alpha_\ell z^\ast-(\alpha_\ell \lbbp + \ov{c}_0)$ and $\alpha_\ell=w_\ell + (\mmax-1)\tveh$.
        Based on Corollary \ref{cor1}, if $|\lcb{R}| < \dmax$, then $\lcb{R}$ contains the routes of any optimal solution and is defined \textit{optimal}.
        Problem $\msp(\lcb{R})$ is obtained from problem $\msp$ where the route set $\lc{R}$ is substituted with $\lcb{R}$.

    \item \textit{Solve problem $\msp(\lcb{R})$}. Let $\ov{\mb{x}}$ be the solution obtained and let  $\ov{z}(\msp(\lcb{R}))$ be its cost (we assume that $\lcb{R}$ contains also the routes corresponding to the current solution $\mb{x}^\ast$); we impose a time limit of \tlim\ seconds in solving $\msp(\lcb{R})$.
        Problem $\msp(\lcb{R})$ is defined \textit{optimal} if it has been solved to optimality within the imposed time limit; otherwise it is defined \textit{not optimal} (see Section \ref{sec:solfr}). If $z^\ast > \ov{z}(\msp(\lcb{R}))$ set $z^\ast = \ov{z}(\msp(\lcb{R}))$ and $\mb{x}^\ast=\ov{\mb{x}}$.
    \item \textit{Test if the $\msp(\lcb{R})$ solution obtained is an optimal \vrpfof\ solution}.
        Let $\gmin$ be a \rev{dual bound} on the reduced cost of any route that has not been generated, i.e., $\ov{c}_\ell \geq \gmin$, $\forall \ell \in \lc{R} \setminus \lcb{R}$,  and let $\lbnew$ be a \rev{dual bound} on the cost of any solution to $\msp$ involving one or more routes of the set $\lc{R} \setminus \lcb{R}$, computed as $\lbnew = \lbbp + (\gmin + \ov{c}_0)/(\mmax \tveh)$. If one of the following two conditions applies, then $\mb{x}^\ast$ is guarantee to be an optimal \vrpfof\ solution and the algorithm terminates:
        \begin{enumerate}[(i)]
          \item Both $\lcb{R}$ and $\msp(\lcb{R})$ are \textit{optimal}, or
          \item $\msp(\lcb{R})$ is \textit{optimal} and $z^\ast \leq \lbnew$.
        \end{enumerate}
    \item \textit{Termination condition}. If $i \geq \itmax$ or $(z^\ast-\lbnew)/\lbnew \leq \gmax$, then Stop.
    \item \textit{Updating $\dmax$ and \tlim}. Set $\dmax=\varepsilon_1 \dmax$ ($\varepsilon_1 > 1$) and $\tlim=\tlim + \varepsilon_2$ ($\varepsilon_2 > 0$) and go to Step 4 ($\varepsilon_1$ and $\varepsilon_2$ are two user-defined parameters).
\end{enumerate}

\er{The exact method starts with the use of both bounding procedures \cb~and \ccg~to compute valid primal and dual bounds for our problem (Steps 2 and 3). As explained in Section~\ref{sec:dual_ascent}, procedure \ccg~provides a good approximation of problem $DNC\msp$. Thanks to Theorem~\ref{thm2} it is hence possible to obtain $(\gb{\mu}, \omega)$, the corresponding dual variables of $DCC\msp$.
These new duals allow to compute the reduces costs $\ov{c}_\ell$, $\forall \ell \in \lc{R}$.\\
Step 4 selects the $\dmax$ columns to be used in the solution of $\msp(\lcb{R})$, Corollary~\ref{cor1} allows to restrict set $\lcb{R}$ to only useful columns. It is worth noting that the primal bounds computed in Steps 1 and 2 play an important role in the identification of such columns.
In Step 5, problem $\msp(\lcb{R})$ is solved to optimality (up to a given time limit). The exact procedure used to solve it is explained in Section~\ref{sec:solfr}. If the primal solution is better than the best primal solution found so far, $z^*$ is updated.\\
In Step 6, the optimality of the solution of $\msp(\lcb{R})$ is checked. If the number of columns added to $\lcb{R}$ in Step 5 was lower than $\dmax$ (i.e., $\lcb{R}$ is optimal) and we were able to solve $\msp(\lcb{R})$ within the given time limit, then the algorithm terminates.
On the other hand, if $\msp(\lcb{R})$ is solved to optimality but some of the columns with $\ov{c}_\ell < \gamma_\ell$ has been excluded from $\lcb{R}$, a new dual bound $\lbnew$ is computed, based on the columns that have not been added to $\lcb{R}$. If $z^\ast \leq \lbnew$, the optimal solution of $\msp(\lcb{R})$ is also an optimal solution of $\msp(\lc{R})$.\\
Step 7 terminates the algorithm  if  either the number of maximum iterations is reached or the relative gap between the best primal and dual bounds is sufficiently small.
Finally, in Step 8, $\dmax$ and the time limit are updated.}

Notice that whenever the algorithm terminates at Step 7, problem $\msp$ has not been solved to optimality and value $z^\ast$ represents the cost of the best solution found. The next section describes the method used to solve problem $\msp(\lcb{R})$ whereas the procedure used to generate the reduced set of routes $\lcb{R}$ is described in Section \ref{sec:genr}.

\subsection{Solving problem $\msp(\lcb{R})$ to optimality}\label{sec:solfr}

Problem $\msp(\lcb{R})$ can be rewritten as the following mixed integer linear programming problem  using the transformation proposed by \citet{CC62}:
 \begin{alignat}{3}
(\msp(\lcb{R}))\quad  z(\msp(\lcb{R}))= \min  & \sum_{ \ell \in \lcb{R}} c_\ell y_\ell \label{FR.objective}\\
                    s.t. & \eqref{CF.degree-f}-\eqref{CF.yu},\\
                    & y_\ell \leq u, \quad & \forall \ell \in \lcb{R} \label{FR.y}\\
                    & u - M (1-x_\ell) \leq y_\ell \leq M x_\ell, \quad & \forall \ell \in \lcb{R} \label{FR.yx}\\
                         & u \geq 0, \label{FR.u} \\
                         & y_\ell \geq 0, x_\ell \in \{0,1\} \quad &\forall \ell \in \lcb{R}                       \label{FR.x},
\end{alignat}
where in constraints \eqref{CF.degree-f}-\eqref{CF.yu} set $\lc{R}$ is substituted with $\lcb{R}$ and $M$ is a very large positive number.
The above formulation is impractical to solve, even for moderate size \vrpfof\ instances, since the numbers of variables and of the additional constraints \eqref{FR.y} and \eqref{FR.yx} can be huge.

We now describe an exact method for solving problem $\msp(\lcb{R})$ based on the extension of the method proposed by \citet{Dinkelbach1967} to integer \rev{problems}.
For sake of notation, we rewrite problem $\msp(\lcb{R})$ as
  $(\msp(\lcb{R}))\quad  z(\msp(\lcb{R}))= \min \mml n(\mb{x})/d(\mb{x}): \mb{x} \in P \mmr$,
where $n(\mb{x})=\sum_{\ell \in \lcb{R}}c_\ell x_\ell$, $d(\mb{x})=\sum_{\ell \in \lcb{R}}w_\ell x_\ell$,
$X=\{\mb{x} \in \Real_+^{|\lcb{R}|}: \eqref{F.degree-f}, \eqref{F.degree-c}, \text{ and } \mmin \leq \sum_{\ell \in \lcb{R}}x_\ell \leq \mmax\}$ and
$P=X \cap \{0,1\}^{|\lcb{R}|}$.
We assume that $d(\mb{x}) > 0$, $\forall \mb{x} \in P$, and that $\msp(\lcb{R})$ admits a finite optimal solution.

\revb{The exact method is based on the following theorem \citep[see][]{Dinkelbach1967}.}
\begin{thm}\label{thm4}
$\ov{r}=n(\ov{\mb{x}})/d(\ov{\mb{x}})=z(\msp(\lcb{R}))$ if, and only if, for the parametric problem $FP(r)$,  $z(FP(r))=\min \mml n(\mb{x})-r d(\mb{x}): \mb{x} \in P \mmr$,
$z(FP(\ov{r}))=n(\ov{\mb{x}})-\ov{r}d(\ov{\mb{x}})=0$.
\end{thm}
\noindent {\textit{Proof.}} The proof is provided in the e-companion to this paper. $\Box$

The scheme of the proposed exact method for solving $\msp(\lcb{R})$ is as follows.
\begin{enumerate}[Step 1.]
    \item \textit{Initialization}. Set $\mb{x}^0=\mb{x}^\ast$, where  $\mb{x}^\ast$ is the current best know solution (see Step 5 of the exact method). Set $i=0$.
    \item \textit{Compute the current ratio $r_{i+1}$}. Set $r_{i+1}=n(\mb{x}^i)/d(\mb{x}^i)$ and set $i=i+1$.
    \item \textit{Solve the parametric problem}. Solve problem
             $FP(r_i)$,  $z_i(FP(r_i))=\min \mml n(\mb{x})-r_i d(\mb{x}): \mb{x} \in P \mmr$,
            and let $\mb{x}^i$ be the solution obtained.
    \item \textit{Termination condition}.
            If $z_i(FP(r_i)) < 0$, go to Step 2; otherwise,          set $\ov{z}(\msp(\lcb{R}))=r_i$, $\ov{\mb{x}}=\mb{x}^i$, $\ov{r}=r_i$ and if solution $\mb{x}^i$ is an optimal solution of problem $FP(r_i)$, define
            $\msp(\lcb{R})$ \textit{optimal}; otherwise define
                  $\msp(\lcb{R})$ \textit{not optimal}.
\end{enumerate}

In solving problem $FP(r_i)$, the time limit of \tlim\ seconds introduced at Step 5 of the exact method is imposed, therefore solution $\mb{x}^i$  is a proven optimal solution of $FP(r_i)$ if the problem has been solved within the imposed time limit.

\revb{The following Lemma \ref{lemmaexm} revisits a result from \citet{Espinoza2010} to show the correctness of the above iterative algorithm by also considering the termination conditions of Step 4.}

\begin{lem}\label{lemmaexm} The following properties hold about the exact algorithm for solving $\msp(\lcb{R})$.
  \begin{enumerate}[(1)]
    \item The sequence $\{r_i\}$ is monotone decreasing, i.e., $r_i > r_{i+1}$, for all $i$ such that $z_i(FP(r_i)) <0$.
    \item \rev{If $z_{i}(FP(r_i))\geq0$ and solution $\mb{x}^i$ is optimal}, the value $\ov{z}(\msp(\lcb{R}))$ corresponds to the optimal solution value and $\ov{r}=n(\ov{\mb{x}})/d(\ov{\mb{x}})$.
    \item \rev{If $z_{i}(FP(r_i))<0$, we have} $d(\mb{x}^i) > d(\mb{x}^{i+1})$.
  \end{enumerate}
\end{lem}
\noindent {\textit{Proof.}} The proof is provided in the e-companion to this paper. $\Box$

\revb{Based on the lemma above, the following theorem
shows that the convergence rate of the algorithm is superlinear \citep{Espinoza2010}.}

\begin{thm}\label{thm3}
  For all $r_i \neq \ov{r}$ we have
  \begin{equation}\label{thm3-1}
    \frac{\ov{r}-r_{i+1}}{\ov{r}-r_i} \leq 1 - \frac{d(\ov{\mb{x}})}{d(\mb{x}^i)} < 1.
  \end{equation}
\end{thm}
\noindent {\textit{Proof.}} The proof is provided in the e-companion to this paper. $\Box$

\section{Pricing problem and generation of sets $\lcb{R}$}\label{sec:ng-genr}

In this section, we describe the details of the pricing problem associated with bounding procedures  \dk, \cg\ and \ccg\ and the procedure used to compute functions $\varphi^t_i$ in bounding procedure \cb. Moreover, we describe the details of the procedure used to generate sets $\lcb{R}$ in the exact method.

\subsection{Route relaxation $ng$-routes}\label{sec:ngroute}

\nomenclature[ch]{$N_i$, $\deltap{N_i}$, $\Pi(P)$, $\lc{E}$, $\Psi$, $\ov{d}_{ij}$ }{ng-route}

The pricing problem associated with procedure \dk, \cg\ and \ccg\ requires to find minimum cost elementary routes over a graph with both positive and negative edge and arc costs, a strongly $\mathcal{NP}$--hard problem. Therefore, in practice we enlarge the set of routes $\lc{R}$ to contain also \textit{non-necessarily elementary} routes, i.e., coefficients $a_{i\ell}$ are general nonnegative integers. Although non-elementary routes are infeasible, this relaxation has the advantage that the pricing subproblem becomes solvable efficiently (by \dynp). Moreover, Theorem \ref{thm1} remains valid if the set of routes $\lc{R}$ is enlarged to contain also non-necessarily elementary routes. The relaxation we used is based on the route relaxation proposed by \citet{Baldacci2011} for the \vrptw\ and can be described as follows.

Let $N_i \subseteq \setfc$ be a set of selected customers for vertex $i$ (according to some criterion) such that $N_i \ni i$ and $\vert N_i \vert \leq \deltap{N_i}$, where $\deltap{N_i}$ is a parameter (e.g., $\deltap{N_i}=5$, $\forall i \in \setfc$, and $N_i$ contains $i$ and the four nearest customers to $i$). The sets $N_i$ allow us to associate with each forward path $P=(0,i_1,\dots,i_k)$ the subset $\Pi(P) \subseteq V(P)$, $V(P)=\{0,i_1,\ldots,i_{k-1},i_k\}$, containing customer $i_k$ and every customer  $i_r$, $r=1,..,k-1$, of $P$ that belongs to all sets  $N_{i_{r+1}},\dots,N_{i_k}$  associated with the customers  $i_{r+1},\dots,i_k$  visited after $i_r$. The set $\Pi(P)$ is defined as:
    $\Pi(P) = \{ i_r: i_r \in \bigcap_{s=r+1}^k N_{i_s}, r=1,\dots,k-1\} \cup \{i_k\}$.
A \textit{$ng$-path} $(NG,t,i)$ is a non-necessarily elementary path $P=(0,i_1,\dots,i_{k-1},i_k=i)$ starting from the depot at time 0, visiting a subset of customers (even more than once) such that $NG=\Pi(P)$, ending at customer $i$ at time $t$, and such that $i \notin \Pi(P')$, where $P'=(0,i_1,\dots,i_{k-1})$ is an $ng$-path. We denote by $f(NG,t,i)$ the cost of the least cost $ng$-path $(NG,t,i)$.
An $(NG,t,i)$-route is an $(NG,t,0)$-path visiting at time $t$ the last customer $i$ before arriving at the depot. The cost of the least cost $(NG,t,i)$-route is given by $f(NG, t,i) + d_{i0}$.
Functions $f(NG,t,i)$ can be computed using \dynp\ as follows.
The state space graph $\lc{H}=(\lc{E},\Psi)$ is defined as follows:
$\lc{E} = \{(NG,t,i) \, : \, \forall NG \subseteq N_i \text{ s.t. } NG \ni i, \, \forall t, 0 \leq t \leq \tveh, \forall i \in \setv \}$,
$\Psi=\{((NG',t',j),(NG,t,i)): \forall (NG',t',j) \in \Psi^{-1}(NG,t,i), \forall (NG,t,i) \in \lc{E} \}$,
where $\Psi^{-1}(NG,t,i)=\{(NG',t-s_i-t_{ji},j): \forall NG' \subseteq N_j \text{ s.t. }NG' \ni j \text{ and } NG' \cap N_i=NG \setminus\{i\}, \,t-s_i-t_{ji} \geq 0, \, \forall j \in \setv \setminus \{i\} \}$.

The \dynp\ recursion for computing $f(NG,t,i)$ is as follows:
\begin{equation}
f(NG,t,i) = \min_{(NG',t',j) \in \Psi^{-1}(NG,t,i)}\{f(NG',t',j)+d_{ji}\}, \quad \forall (NG,t,i) \in \lc{E}. \label{ng.a11}
\end{equation}
The following initialization is required: $f(\{0\},0,0)=0$ and $f(\{0\},t,0)= \infty$,  $\forall t$ such that $0 < t \leq \tveh$.



%
%

\subsubsection*{Computing functions $\varphi^t_i$ at Step 2 of algorithm \cb}

We first compute functions $f(NG,t,i)$ using \dynp\ recursion \eqref{ng.a11} and the modified costs $\ov{d}_{ij}$ instead of $d_{ij}$, where $\ov{d}_{ij}=d_{ij}-\lambda_j$. Then we compute functions $\varphi^t_i$, $\forall i \in \setfc$, $\forall t, \max \mml 1,\tmin-(\nveh-1)\tveh \mmr \leq t \leq \tveh$, as
          $\varphi^t_i = \min_{(NG,t-t_{i0},i) \in \lc{E}} \mml f(NG,t-t_{i0},i) + \ov{d}_{i0}\mmr$.

\subsection{Generating set $\lcb{R}$: procedure \genr}\label{sec:genr}

The generation of the reduced route set $\lcb{R}$ performed at Step 4 is based on a similar procedure proposed by \citet{Baldacci2011} for the \vrp\ with \rev{Time Windows}, used to generate elementary and feasible routes.
Given a dual solution $(\gb{\mu}, \omega)$ of $DCC\msp$  and a user defined parameter $\dmax$, the procedure generates the largest subset $\lcb{R} \subseteq \lc{R}$ satisfying conditions \eqref{DynP.Conds}-a and \eqref{DynP.Conds}-b. Procedure \genr\ is a \dynp\ programming that is analogous to \rev{Dijkstra's} algorithm on an expanded state-space graph dynamically generated.

Associate with each arc $(i, j)\in A$ modified arc cost $\ov{d}_{ij}$ defined as $\ov{d}_{ij} = d_{ij} - \mu_j - (t_{ij}+s_j)\omega$.
It is easy to see that the reduced cost with respect to the dual vector $\gb{\mu}$ and
value $\omega$ can be computed as $\ov{c}_\ell = \sum_{(i,j) \in A(R_\ell)} \ov{d}_{ij}$. Procedure \genr\ dynamically generates a
state-space graph where each state corresponds to a feasible \textit{forward path}.

A \textit{forward path} $P=(0,i_1,\ldots,i_{k-1},i_k)$ is an elementary path starting from depot 0 at time $0$, visiting vertices $V(P)=\{0,i_1,\ldots,i_{k-1},i_k\}$ and ending at customer  $i_k = \ep{P}$ at time $t(P)$ with $t(P) \leq \tveh$. We denote by $A(P)$ the set of arcs traversed by path $P$ and by $c(P)= \sum_{(i,j)\in A(P)}d_{ij}$ the cost of path $P$.
Let $DB(P)$ be a \rev{dual bound} on the reduced cost of any route that contains a forward path $P$. Any forward path $P$ such that $DB(P) \geq \gamma$ cannot be part of a route in the set $\lcb{R}$ that satisfies conditions \eqref{DynP.Conds}, where $\gamma=\alpha z^\ast-(\alpha \lbbp + \ov{c}_0)$ and $\alpha= \mmax\tveh$.

Let $\tau$ be a set of temporary feasible forward paths that is initialized by setting $\tau = \{P_0\}$, where $P_0$ represents the initial empty path such that $\ep{P_0} = 0$ and $t(P_0)=0$. The route set $\lcb{R}$ is initialized by setting $\lcb{R}=\emptyset$.
At each iteration of algorithm \genr\ the forward path $P \in \tau$ having the smallest \rev{dual bound} value (i.e., such that $DB(P) = \min \mml DB(P): P \in \tau \mmr$ ) is extracted from the set $\tau$. The expansion of a forward path $P$ are derived by extending $P$ with arc $(\ep{P},j) \in A$, $\forall j \notin V(P)\setminus \{0\}$. We have two cases:
\begin{enumerate}[i)]
  \item $j=0$. The expansion of forward path $P$ creates a route; if the route is feasible and satisfy condition \eqref{DynP.Conds}-b it is inserted in the set $\lcb{R}$;
  \item $j \neq 0$. The expansion of path forward path $P$ creates a path $P'$; if the path is feasible and $DB(P') < \gamma$ (see above) it is added to the set $\tau$.
\end{enumerate}

Procedure \genr\ terminates when either $\tau = \emptyset$  or $|\lcb{R}| = \dmax$. Since the size of the set $\tau$ is exponential, we impose that the size of the set $\tau$ cannot exceed an
a-priori defined limit $\nstatb$. If $|\tau|$ becomes greater than $\nstatb$, procedure \genr\ terminates prematurely. At the end of the procedure, value \gmin\ is set equal to $\max_{P \in \tau}\mml DB(P) \mmr$ if $|\lcb{R}| = \dmax$ or $|\tau|=\nstatb$, and $\infty$ otherwise.

\subsubsection{Computing $DB(P)$}\label{sec:genr_lb}

We define a \textit{backward path} $\ov{P}=(\ep{\ov{P}}=i_k, i_{k+1}, \dots, i_h, 0)$ as a path starting from vertex $\ep{\ov{P}}$ at time $t(\ov{P})$, visiting vertices in $V(\ov{P})=\{i_k, i_{k+1}, \dots, i_h, 0\}$ and ending at the depot before time $\tveh$.

A \rev{dual bound}s on the cost $c(\ov{P})$ of $\ov{P}$ can be computed using the $ng$-path by defining nonnecessarily elementary \textit{backward} $ng-path$ $(NG,t,i)$ similarly to the forward \textit{$ng$-path} $(NG,t,i)$ defined in Section \ref{sec:ngroute}. Let $f^{-1}(NG,t,i)$ be the cost of the least-cost backward $ng-path$ $(NG,t,i)$. Functions $f^{-1}(NG,t,i)$ can be computed with the same \dynp\ recursions used to compute $f(NG,t,i)$ by replacing the cost and time matrices $[d_{ij}]$ and $[t_{ij}]$ with their transposed matrices $[d_{ij}]^{T}$ and $[t_{ij}]^{T}$.
We have that the cost $c(\ov{P})$ of any elementary backward path $\ov{P}$ satisfies the following inequality:
\begin{equation}
c(\ov{P}) \geq \min_{NG \subseteq V(\ov{P}) \cap N_{\ep{\ov{P}}}}\{f^{-1}(NG,t(\ov{P}),\ep{\ov{P}})\}.
\end{equation}

To compute \rev{dual bound} $DB(P)$, function $f^{-1}(NG,t,i)$ are computed with the modified costs $[\ov{d}_{ij}]$, and the subsets $N_i$, $i \in \setfc$, contain the $\deltap{N_i}$ nearest customers to $i$ according to $\ov{d}_{ij}$. \rev{Dual bound} $DB(P)$ is computed as follows:
\begin{equation}
 DB(P) = \sum_{(i,j)\in A(P)}\ov{d}_{ij} + \min_{\substack{NG \subseteq N_{\ep{P}} \text{ s.t. } NG \cap V(P) = \{\ep{P}\} \\ t' \leq \tveh-t(P)}}  \{ f^{-1}(NG,t',\ep{P}) \}.
\end{equation}

\subsubsection{Dominance rules}\label{sec:genr_dom}

A speed-up in procedure \genr\ can be obtained by removing dominated
paths from the set $\tau$. A dominated path is either a path that cannot lead to a feasible route or a path such that any route containing it cannot be part of any optimal solution.
Dominance rules are defined based on the type of fractional objective function considered as follows.

\begin{dom}{Minimization of Cost/Load}
  A forward path $P_1$ dominates a forward path $P_2$ if $\ep{P_1}=\ep{P_2}$, $V(P_1)=V(P_2)$ and $c(P_1) \leq c(P_2)$.
\end{dom}

\begin{dom}{Maximization of Profit/Time}
  A forward path $P_1$ dominates a forward path $P_2$ if $\ep{P_1}=\ep{P_2}$, $V(P_1)=V(P_2)$ and $t(P_1) \leq t(P_2)$.
\end{dom}

\section{Computational results}\label{sec:computational results}

This section reports on the computational results of the dual and primal bounds and the exact method described in this paper.
All algorithms were coded in C++ and compiled with Microsoft Visual Studio 2013 compiler. The IBM ILOG CPLEX 12.6.4 callable library (\citet{CPLEX}) was used as the integer programming solver for solving the parametric problem $FP(r_i)$ in the exact method (see Section \ref{sec:solfr}).
All tests were run on a Lenovo ThinkStation P300 (i7-4790 CPU @ 3.6 GHz - 32 GB of RAM) running under Microsoft Windows 7 Professional operating system.

\subsection{\rev{Instances description}}
\er{The bounding procedures and the exact method were tested on two classes of instances, namely \classCA and \classPA, corresponding to \vrpfof\ instances with the two objective functions (described in Section \ref{sec:problem description}) minimization of Cost/Load and maximization of Profit/Time respectively.}

The instances of the \er{two} classes were derived from instances proposed in the literature for the Capacitated \vrp\ (\cvrp). More precisely, we considered the 27 instances, and corresponding optimal solutions, of class A generated by \citet{A95} and available at \url{http://vrp.galgos.inf.puc-rio.br/index.php/en/}.

Let \cvrp($\ov{n}$,$[q_i]$,$K$,$Q$,$[c_{ij}]$) be a \cvrp\ instance, where $\ov{n}$ represents the number of vertices (including the depot), $[q_i]$ the customer demands, $K$ the number of vehicles, $Q$ the vehicle capacity, and $[c_{ij}]$ the cost matrix; cost matrix $[c_{ij}]$ is computed according to the TSPLIB EUC\_2D standard (see \citet{Reinelt1991}).
For each \cvrp($\ov{n}$,$[q_i]$,$K$,$Q$,$[c_{ij}]$) instance of class A, we generate an instance for each of \er{two} classes as follows.
\begin{enumerate}[i)]
  \item The depot and the set of customers correspond to the depot and the set of customers of the original \cvrp\ instance. We set $\nf=\lfloor \alpha (\ov{n}-1)\rfloor$, $\alpha < 1$;
  \item Class \classCA. We set $d_{ij}=c_{ij}$, $\forall (i,j) \in A$, $t_{ij}=0$, $\forall (i,j) \in A$, $\tveh=Q$, $s_i=q_i$, $\forall i \in \setfc$, and $m=\min \mml BPP(\ov{n}-1,[q_i],Q)+1, K\mmr$, where $BPP(\ov{n}-1,[q_i],Q)$ is the cost of the optimal solution of the Bin Packing Problem instance with $\ov{n}-1$ items, weights $[q_i]$ and bin capacity equal to $Q$.
  \item Class \classPA. We set $d_{ij}=-q_i$, $\forall (i,j) \in A$, $t_{ij}=c_{ij}$, $\forall (i,j) \in A$.
\er{Service times $\{s_i\}$, maximum working time $\tveh$ and maximum number of vehicles $\nveh$ are computed using the best solution found for the corresponding \cvrp\ instance and a simple heuristic algorithm for the \vrpfof, used to guarantee the feasibility of the corresponding instance (details are omitted for sake of brevity).}
\end{enumerate}

\begin{table}
\caption{Parameters used by the different procedures} \label{tab:param}
\begin{center}
{\footnotesize{
\setlength{\tabcolsep}{2.7pt}
\renewcommand{\arraystretch}{1.0}
\begin{tabular}{l|l}
\hline
Procedure & Parameters \bigstrut\\
\hline
\ccg  & $\pi_i=s_i, \forall i \in \setfc$ \\
Pricing $(NG,t,i)$-routes  & $\deltap{N_i}=12$ nearest nodes to $i$ according to cost matrix $[d_{ij}]$ \\
\cb   & $\epsilon=1.0$,  $Maxit3=200$ \\
\genr & $\nstatb=200E+6$ \\
Exact method & $\itmax=3$, $\dmax=300,000$, $\tlim=3,600$, $\gmax=\infty$, $\varepsilon_1=5$, $\varepsilon_2=3,600$ \bigstrut[b]\\
\hline
\end{tabular}%
}}
\end{center}
\end{table}

In generating the instances, we used $\alpha \in \{0.5,0.75\}$, therefore two instances per class are generated for each \cvrp\ instance. Given the original \cvrp\ instance name ``$<name>$'', the instance with $\alpha=0.5$ is denoted with $<name>a$ whereas the instance with $\alpha=0.75$ is denoted with $<name>b$.
A total number of 162 instances were generated, 54 instances per class. All the instances are available upon request to the authors as text files.

Based on the results of preliminary experiments to identify good parameter settings for our algorithms, we decided to use the settings reported in Table \ref{tab:param}. The following section reports on the results about the dual and primal bounds computed by procedures $\ccg$ and $\cb$ whereas Section \ref{sec:cr_lb} shows the results obtained by the exact method.

\subsection{Computational results on the dual and primal bounds}\label{sec:cr_lb}

To compare dual bounds $\lbcb$ and $\lbbp$ computed by procedures \cb\ and and \ccg, respectively, we implemented a standard column generation algorithm (called \cg) based on formulation $CC\msp$  and Dinkelbach's algorithm (called \dk) applied to formulation $C\msp$.
Both algorithms \cg\ and \dk\ are based on the $(NG,t,i)$-routes relaxation and the linear programming solver of IBM CPLEX is used to solve the master problem at each iteration of algorithm \cg.


Table \ref{tab:sumlb} summarises the results obtained about the different dual and primal bounds on the \er{two} classes of instances. For bounding procedures $\cb$, $\ccg$ and \cg, the table reports the average percentage deviation of the dual bound (column $\%B$) and the average computing time in seconds (column ``Time''). The percentage deviation is computed as $100.0 \er{\times} B/z^\ast$, where $z^\ast$ is the cost of the best solution found and $B$ is the value of the dual bound;
for class \classCA, $B$ refers to a \er{lower bound} whereas for class \classPA, $B$ refers to a \er{upper bound}. For bounding procedure \dk, the table reports only the average computing time in seconds, being the value of the dual bound computed by the procedure equal to the one computed by procedure \cg.

For bounding procedures $\cb$ and $\ccg$, the table also reports the average percentage deviation of the primal bound obtained by the heuristic procedure (column $\%PB$), computed as $100.0 \er{\times} PB/z^\ast$, where $z^\ast$ is the cost of the best solution found and $PB$ is the value of the primal bound; column ``Time'' also includes the time spent for computing the primal bound and, in the case of procedures $\cb$, the time spent for computing value $\tmin$.

For procedures \ccg, \cg, and \dk, the table also shows the average number of columns or variables of the final master problem. In particular, for procedure \dk\ the number is computed as the sum over all algorithm iterations, whose average number is reported under column ``Iter'' in the table.

Complete computational results about the dual and primal bounds are reported in the e-companion to this paper. Moreover, the e-companion also reports statistics about the best solutions found by the heuristic procedures and the exact method, including the values $z^\ast$ used to compute the different percentage deviations.

\begin{table}
\caption{Summary results on the dual bounds} \label{tab:sumlb}
\begin{center}
{\footnotesize{
\setlength{\tabcolsep}{2.0pt}
\renewcommand{\arraystretch}{1.0}
\begin{tabular}{r|rrr|rrrr|rrr|rrr}
\hline
Class & \multicolumn{3}{c|}{Procedure \cb} & \multicolumn{4}{c|}{Procedure \ccg}     & \multicolumn{3}{c|}{Procedure \cg} & \multicolumn{3}{c}{Procedure \dk} \bigstrut[t]\\
      & \multicolumn{1}{c}{$\%B$} & \multicolumn{1}{c}{$\%PB$} & \multicolumn{1}{c|}{\tabtime} & \multicolumn{1}{c}{$\%B$} & \multicolumn{1}{c}{$\%PB$} & \multicolumn{1}{c}{\tablbcols} & \multicolumn{1}{c|}{\tabtime} & \multicolumn{1}{c}{$\%B$} & \multicolumn{1}{c}{\tablbcols} & \multicolumn{1}{c|}{\tabtime} & \multicolumn{1}{c}{\tablbit} & \multicolumn{1}{c}{\tablbcols} & \multicolumn{1}{c}{\tabtime} \bigstrut[b]\\
\hline
\classCA & 98.5  & 107.0 & 18.8  & 98.7  & 106.3 & 374.2 & 0.4   & 98.8  & 7683.7 & 1.2   & 2.8   & 7254.2 & 1.1 \\
\classPA & 106.8 & 96.3  & 88.9  & 106.1 & 97.7  & 610.1 & 0.7   & 105.9 & 2282.0 & 0.7   & 2.6   & 2377.5 & 0.8 \bigstrut[b]\\
\hline
\end{tabular}%
}}
\end{center}
\end{table}


Table \ref{tab:sumlb} shows that the dual bounds computed are on average quite tight for \er{instances belonging  to the} class \classCA. Clearly, the dual bound computed by procedure \cg\ dominates the dual bound produced by procedure \ccg\ and it is equivalent to the one computed by procedure \dk. Instances of class \classPA\ are more difficult for our bounding procedures, as shown  by the average percentage deviation reported in the table; this is probably due to cost structure of the class.
The dual bounds computed by procedure \ccg\ are on average very close to the ones produced by \cg\ and can be computed faster. Column $|\lcb{R}|$ (i.e, the average number of columns of the final master problem) shows that procedure \ccg\ is not affected by the typical degeneracy of standard column generation generation based on the simplex, like \cg. The detailed results reported in the e-companion show that the dual bound produced by \cb\ is always inferior with respect to one produced by \cg. Further, its computation is more time consuming due to the complexity of the corresponding \dynp\ recursion - it is worth mentioning that the computation of value $\tmin$ is negligible.

Concerning the primal bounds computed by procedures \cb\ and \ccg, the table show that both the two procedures can compute good quality solutions.
The detailed results reported in the e-companion also show that it is convenient to compute both primal bounds. The detailed results show that the heuristic algorithm applied during the execution of procedure  \ccg\ failed to compute feasible primal bounds for three instances of class \classPA.

\subsection{Computational results on the exact method}\label{sec:cr_lb}

Tables \ref{tab:exm-CA}-\ref{tab:exm-PA} show the results about the exact method.
For each instance, we report a symbol to denote if the instance was solved to optimality (``\tabopt''), and the corresponding total computing time in seconds, that also includes the time spent by the bounding procedures executed at steps 2 and 3 of the exact method (see Section \ref{sec:exm}) and the time spent by procedure \genr.

The next three blocks of columns show the details of the iterations of the exact method. For each iteration, we report the
cardinality of the reduced set $\lcb{R}$ ($|\lcb{R}|$) generated by procedure \genr, the percentage deviation of the value of the optimal solution of the reduced integer program $\msp(\lcb{R})$ ($\%z^\ast$), a symbol (``\tabipnotopt'') to denote if the time limit imposed was reached in solving $\msp(\lcb{R})$ (IP), the number of iterations executed to solve problem  $\msp(\lcb{R})$ with the procedure described in Section \ref{sec:solfr}, the percentage deviation of bound $\lbnew$ with respect to $\%z^\ast$ ($\%B$) (see Step 6 of the exact method) and the computing time in seconds spent to solve problem $\msp(\lcb{R})$.

In the tables, under column $|\lcb{R}|$, symbols ``\tabnstatb'' and ``\tabmaxcols'' are reported whenever limit $\nstatb$ or $\dmax$ of procedure \genr\ has been reached, respectively. Further, the heading of each table shows, for each iteration of the method, the value of parameters  $\dmax$ and \tlim.

Details about the best solution found by the heuristic procedures and the exact method are reported in the e-companion. In particular, for each instances, it is shown the values $z^\ast$ of the best solution found (including the corresponding numerator and denominator), the number of routes of the solution, the number of optional customers selected in solution, the percentage of the working time utilization and the average number of customers per route.

Tables \ref{tab:exm-CA}-\ref{tab:exm-PA} show that 53 and 30 out of the 54 instances per class were solved to optimality by the exact method for classes \classCA\ and \classPA, respectively. For classes \classCA, instances with up to 79 customers were solved to optimality whereas the largest instance solved to optimality for class \classPA\ involves 62 customers. Instances of class \classPA\ are clearly more difficult for our method, as also shown by the quality of the corresponding dual bounds, and by the fact that the limit $\dmax$ imposed at the different iterations to the exact method is generally reached. For the instances not solved to optimality for classes \classCA, the final bound $\%B$ is very tight, thus showing that near-optimal solutions are also computed for the corresponding instances.
Notice that in the tables, whenever under column $|\lcb{R}|$ of the last iteration appears symbol ``\tabmaxcols'' (i.e., $\dmax$ limit reached) and the instance has been solved to optimality, the condition on the dual bound $\lbnew$ (see Step 6 of the exact method) is used as optimality condition.
Most of the instances of class \classCA\ can be solved to optimality within the first two iterations of the exact method and the limit generally attained during the different iterations is the maximum number of columns or routes $\dmax$. In particular,
46 and 12 instances were solved to optimality at the first iteration of the exact algorithm for the \er{two} classes of instances, respectively, thus showing the effectiveness of our iterative exact procedure in reducing the size of the integer problems solved at each iteration of the Dinkelbach's algorithm.

In order to have some insights about the type of instances used and the solutions computed, the e-companion reports some statistics about the best solutions found for the \er{two} classes of instances. In particular, the tables show that the maximum number of vehicles $\nveh$ is generally tight and that the percentage of the working time utilization is on average superior to 90\%; optional customers are generally included in the best solutions found.



\begin{landscape}
\renewcommand{\arraystretch}{0.7}
{\tiny{
\begin{longtable}{r|r r|r r r r r r|r r r r r r |r r r r r r}
\caption{Instances of class \classCA: exact method} \label{tab:exm-CA}\\
\hline
\multicolumn{3}{c|}{} &  \multicolumn{6}{c|}{$Iter=1$, $\dmax=300,000$, $\tlim=3,600$} & \multicolumn{6}{c|}{$Iter=2$, $\dmax=1,500,000$, $\tlim=7,200$}   &  \multicolumn{6}{c}{$Iter=3$, $\dmax=7,500,000$, $\tlim=10,800$} \bigstrut\\
\hline
\tabname  & Opt   & Time  & $|\lcb{R}|$     & $\%z^\ast$   & \taboptip    & Iter  & \tablb & Time  & $|\lcb{R}|$     & $\%z^\ast$   & \taboptip    & Iter  & \tablb & Time  & $|\lcb{R}|$     & $\%z^\ast$   & \taboptip    & Iter  & \tablb & Time \bigstrut\\
\endfirsthead
\multicolumn{3}{@{}l}{\ldots continued}\\
\hline
\multicolumn{3}{c|}{} &  \multicolumn{6}{c|}{$Iter=1$} & \multicolumn{6}{c|}{$Iter=2$}   &  \multicolumn{6}{c}{$Iter=3$} \bigstrut\\
\hline
\tabname  & Opt   & Time  & $|\lcb{R}|$     & $\%z^\ast$   & \taboptip    & Iter  & \tablb & Time  & $|\lcb{R}|$     & $\%z^\ast$   & \taboptip   & Iter  & \tablb & Time  & $|\lcb{R}|$     & $\%z^\ast$   & \taboptip    & Iter  & \tablb & Time \bigstrut\\
\hline
\endhead \\
\multicolumn{21}{c}{\tabopt: solved to optimality \hspace{0.5cm} \tabnstatb: $\nstatb$ limit reached \hspace{0.5cm} \tabmaxcols: $\dmax$ limit reached \hspace{0.5cm} \tabipnotopt: \tlim\ limit reached }\bigstrut\\
\endfoot
\hline
A-n32-k5a & \tabopt & 8.4   & 24,022 & 100.0 &       & 2     &       & 1.4   &       &       &       &       &       &       &       &       &       &       &       &  \\
A-n32-k5b & \tabopt & 21.0  & 35,174 & 100.0 &       & 2     &       & 2.9   &       &       &       &       &       &       &       &       &       &       &       &  \\
A-n33-k5a & \tabopt & 25.2  & \tabmaxcols & 100.0 &       & 3     &       & 7.7   &       &       &       &       &       &       &       &       &       &       &       &  \\
A-n33-k5b & \tabopt & 20.7  & 632   & 100.0 &       & 2     &       & 0.0   &       &       &       &       &       &       &       &       &       &       &       &  \\
A-n33-k6a & \tabopt & 12.4  & \tabmaxcols & 100.0 &       & 3     &       & 6.0   &       &       &       &       &       &       &       &       &       &       &       &  \\
A-n33-k6b & \tabopt & 24.2  & \tabmaxcols & 100.0 &       & 3     &       & 11.2  &       &       &       &       &       &       &       &       &       &       &       &  \\
A-n34-k5a & \tabopt & 28.1  & \tabmaxcols & 100.0 &       & 2     &       & 18.5  &       &       &       &       &       &       &       &       &       &       &       &  \\
A-n34-k5b & \tabopt & 37.5  & \tabmaxcols & 100.0 &       & 2     &       & 21.3  &       &       &       &       &       &       &       &       &       &       &       &  \\
A-n36-k5a & \tabopt & 38.4  & \tabmaxcols & 100.0 &       & 3     &       & 11.1  &       &       &       &       &       &       &       &       &       &       &       &  \\
A-n36-k5b & \tabopt & 39.5  & \tabmaxcols & 100.0 &       & 2     &       & 19.3  &       &       &       &       &       &       &       &       &       &       &       &  \\
A-n37-k5a & \tabopt & 156.8 & 206,203 & 100.0 &       & 2     &       & 21.3  &       &       &       &       &       &       &       &       &       &       &       &  \\
A-n37-k5b & \tabopt & 75.1  & \tabmaxcols & 100.0 &       & 2     &       & 52.0  &       &       &       &       &       &       &       &       &       &       &       &  \\
A-n37-k6a & \tabopt & 12.1  & 50,360 & 100.0 &       & 2     &       & 4.1   &       &       &       &       &       &       &       &       &       &       &       &  \\
A-n37-k6b & \tabopt & 25.3  & \tabmaxcols & 100.0 &       & 3     &       & 14.0  &       &       &       &       &       &       &       &       &       &       &       &  \\
A-n38-k5a & \tabopt & 137.7 & \tabmaxcols & 100.0 &       & 2     &       & 132.4 &       &       &       &       &       &       &       &       &       &       &       &  \\
A-n38-k5b & \tabopt & 24.7  & \tabmaxcols & 100.0 &       & 2     &       & 14.0  &       &       &       &       &       &       &       &       &       &       &       &  \\
A-n39-k5a & \tabopt & 79.4  & \tabmaxcols & 100.0 &       & 2     &       & 41.9  &       &       &       &       &       &       &       &       &       &       &       &  \\
A-n39-k5b & \tabopt & 35.8  & 8,135 & 100.0 &       & 2     &       & 0.4   &       &       &       &       &       &       &       &       &       &       &       &  \\
A-n39-k6a & \tabopt & 20.7  & \tabmaxcols & 100.0 &       & 2     &       & 9.3   &       &       &       &       &       &       &       &       &       &       &       &  \\
A-n39-k6b & \tabopt & 18.7  & 24,621 & 100.0 &       & 2     &       & 2.1   &       &       &       &       &       &       &       &       &       &       &       &  \\
A-n44-k6a & \tabopt & 18.8  & \tabmaxcols & 100.0 &       & 2     &       & 7.2   &       &       &       &       &       &       &       &       &       &       &       &  \\
A-n44-k6b & \tabopt & 159.3 & \tabmaxcols & 100.0 &       & 2     &       & 142.8 &       &       &       &       &       &       &       &       &       &       &       &  \\
A-n45-k6a & \tabopt & 14.9  & 12,138 & 100.0 &       & 2     &       & 0.4   &       &       &       &       &       &       &       &       &       &       &       &  \\
A-n45-k6b & \tabopt & 91.1  & \tabmaxcols & 100.0 &       & 2     &       & 69.6  &       &       &       &       &       &       &       &       &       &       &       &  \\
A-n45-k7a & \tabopt & 40.7  & 262,658 & 100.0 &       & 2     &       & 27.9  &       &       &       &       &       &       &       &       &       &       &       &  \\
A-n45-k7b & \tabopt & 28.8  & \tabmaxcols & 100.0 &       & 3     &       & 12.0  &       &       &       &       &       &       &       &       &       &       &       &  \\
A-n46-k7a & \tabopt & 10.4  & 6,492 & 100.0 &       & 2     &       & 0.4   &       &       &       &       &       &       &       &       &       &       &       &  \\
A-n46-k7b & \tabopt & 120.9 & \tabmaxcols & 100.0 &       & 3     &       & 113.1 &       &       &       &       &       &       &       &       &       &       &       &  \\
A-n48-k7a & \tabopt & 68.1  & \tabmaxcols & 100.0 &       & 2     &       & 54.5  &       &       &       &       &       &       &       &       &       &       &       &  \\
A-n48-k7b & \tabopt & 91.6  & \tabmaxcols & 100.0 &       & 2     &       & 70.7  &       &       &       &       &       &       &       &       &       &       &       &  \\
A-n53-k7a & \tabopt & 29.5  & \tabmaxcols & 100.0 &       & 2     &       & 4.7   &       &       &       &       &       &       &       &       &       &       &       &  \\
A-n53-k7b & \tabopt & 27.1  & 422   & 100.0 &       & 2     &       & 0.0   &       &       &       &       &       &       &       &       &       &       &       &  \\
A-n54-k7a & \tabopt & 26.0  & \tabmaxcols & 100.0 &       & 2     &       & 9.2   &       &       &       &       &       &       &       &       &       &       &       &  \\
A-n54-k7b & \tabopt & 56.8  & \tabmaxcols & 100.0 &       & 2     & 99.9  & 27.9  & 149,554 & 100.0 &       & 1     &       & 8.5   &       &       &       &       &       &  \\
A-n55-k9a & \tabopt & 21.7  & \tabmaxcols & 100.0 &       & 2     &       & 8.8   &       &       &       &       &       &       &       &       &       &       &       &  \\
A-n55-k9b & \tabopt & 30.6  & \tabmaxcols & 100.0 &       & 3     &       & 24.4  &       &       &       &       &       &       &       &       &       &       &       &  \\
A-n60-k9a & \tabopt & 67.2  & \tabmaxcols & 100.0 &       & 2     &       & 40.7  &       &       &       &       &       &       &       &       &       &       &       &  \\
A-n60-k9b & \tabopt & 68.9  & \tabmaxcols & 100.0 &       & 2     &       & 48.2  &       &       &       &       &       &       &       &       &       &       &       &  \\
A-n61-k9a & \tabopt & 34.0  & \tabmaxcols & 100.0 &       & 3     &       & 16.8  &       &       &       &       &       &       &       &       &       &       &       &  \\
A-n61-k9b & \tabopt & 32.6  & \tabmaxcols & 100.0 &       & 2     &       & 10.9  &       &       &       &       &       &       &       &       &       &       &       &  \\
A-n62-k8a & \tabopt & 368.8 & 99,422 & 100.0 &       & 2     &       & 6.6   &       &       &       &       &       &       &       &       &       &       &       &  \\
A-n62-k8b & \tabopt & 948.4 & \tabmaxcols & 100.0 &       & 2     & 99.1  & 138.5 & \tabmaxcols & 100.0 &       & 1     & 99.6  & 245.5 & 5,681,865 & 100   &       & 1     &       & 507.0 \\
A-n63-k10a & \tabopt & 30.4  & \tabmaxcols & 100.0 &       & 2     &       & 9.8   &       &       &       &       &       &       &       &       &       &       &       &  \\
A-n63-k10b & \tabopt & 95.0  & \tabmaxcols & 100.0 &       & 2     &       & 68.6  &       &       &       &       &       &       &       &       &       &       &       &  \\
A-n63-k9a & \tabopt & 61.6  & \tabmaxcols & 100.0 &       & 2     & 99.4  & 15.2  & 35,241 & 100.0 &       & 1     &       & 0.4   &       &       &       &       &       &  \\
A-n63-k9b & \tabopt & 47.3  & \tabmaxcols & 100.0 &       & 2     &       & 15.5  &       &       &       &       &       &       &       &       &       &       &       &  \\
A-n64-k9a & \tabopt & 118.2 & \tabnstatb & 100.0 &       & 2     & 99.6  & 13.7  & 2,736 & 100.0 &       & 1     &       & 0.0   &       &       &       &       &       &  \\
A-n64-k9b & \tabopt & 118.9 & \tabmaxcols & 100.0 &       & 2     &       & 106.8 &       &       &       &       &       &       &       &       &       &       &       &  \\
A-n65-k9a & \tabopt & 29.4  & \tabmaxcols & 100.0 &       & 2     &       & 9.6   &       &       &       &       &       &       &       &       &       &       &       &  \\
A-n65-k9b & \tabopt & 71.3  & \tabmaxcols & 100.0 &       & 2     &       & 41.7  &       &       &       &       &       &       &       &       &       &       &       &  \\
A-n69-k9a & \tabopt & 120.6 & \tabmaxcols & 100.0 &       & 2     & 99.8  & 75.0  & 6,329 & 100.0 &       & 1     &       & 0.0   &       &       &       &       &       &  \\
A-n69-k9b & \tabopt & 141.9 & \tabmaxcols & 100.0 &       & 2     & 99.2  & 37.1  & 565,327 & 100.0 &       & 1     &       & 51.5  &       &       &       &       &       &  \\
A-n80-k10a & \tabopt & 961.2 & \tabnstatb & 101.0 &       & 2     & 99.5  & 37.7  & \tabmaxcols & 100.0 &       & 2     &       & 431.7 &       &       &       &       &       &  \\
A-n80-k10b & \tabnoopt & 18338.3 & \tabmaxcols & 100.0 &       & 2     & 99.0  & 153.7 & \tabmaxcols & 100.0 & \tabipnotopt & 1     & 99.6  & 7204.3 & \tabmaxcols & 100   & \tabipnotopt & 1     & 99.9  & 10859.8 \\
\hline
\end{longtable}%
}
}
\end{landscape}


\begin{landscape}
\renewcommand{\arraystretch}{0.7}
{\tiny{
\begin{longtable}{r|r r|r r r r r r|r r r r r r |r r r r r r}
\caption{Instances of class \classPA: exact method} \label{tab:exm-PA}\\
\hline
\multicolumn{3}{c|}{} &  \multicolumn{6}{c|}{$Iter=1$, $\dmax=300,000$, $\tlim=3,600$} & \multicolumn{6}{c|}{$Iter=2$, $\dmax=1,500,000$, $\tlim=7,200$}   &  \multicolumn{6}{c}{$Iter=3$, $\dmax=7,500,000$, $\tlim=10,800$} \bigstrut\\
\hline
\tabname  & Opt   & Time  & $|\lcb{R}|$     & $\%z^\ast$   & \taboptip    & Iter  & \tablb & Time  & $|\lcb{R}|$     & $\%z^\ast$   & \taboptip    & Iter  & \tablb & Time  & $|\lcb{R}|$     & $\%z^\ast$   & \taboptip    & Iter  & \tablb & Time \bigstrut\\
\endfirsthead
\multicolumn{3}{@{}l}{\ldots continued}\\
\hline
\multicolumn{3}{c|}{} &  \multicolumn{6}{c|}{$Iter=1$} & \multicolumn{6}{c|}{$Iter=2$}   &  \multicolumn{6}{c}{$Iter=3$} \bigstrut\\
\hline
\tabname  & Opt   & Time  & $|\lcb{R}|$     & $\%z^\ast$   & \taboptip    & Iter  & \tablb & Time  & $|\lcb{R}|$     & $\%z^\ast$   & \taboptip   & Iter  & \tablb & Time  & $|\lcb{R}|$     & $\%z^\ast$   & \taboptip    & Iter  & \tablb & Time \bigstrut\\
\hline
\endhead \\
\multicolumn{21}{c}{\tabopt: solved to optimality \hspace{0.5cm} \tabnstatb: $\nstatb$ limit reached \hspace{0.5cm} \tabmaxcols: $\dmax$ limit reached \hspace{0.5cm} \tabipnotopt: \tlim\ limit reached }\bigstrut\\
\endfoot
\hline
A-n32-k5a & \tabopt & 134.6 & 605   & 100.0 &       & 2     &       & 0.0   &       &       &       &       &       &       &       &       &       &       &       &  \\
A-n32-k5b & \tabopt & 53.5  & 176,046 & 100.0 &       & 2     &       & 2.8   &       &       &       &       &       &       &       &       &       &       &       &  \\
A-n33-k5a & \tabopt & 74.0  & 283   & 100.0 &       & 2     &       & 0.0   &       &       &       &       &       &       &       &       &       &       &       &  \\
A-n33-k5b & \tabopt & 77.3  & 2,285 & 100.0 &       & 1     &       & 0.0   &       &       &       &       &       &       &       &       &       &       &       &  \\
A-n33-k6a & \tabopt & 550.0 & \tabmaxcols & 100.0 &       & 2     & 105.4 & 33.7  & \tabmaxcols & 100.0 &       & 1     & 102.7 & 69.3  & 4,433,401 & 100.0 &       & 1     &       & 387.1 \\
A-n33-k6b & \tabopt & 455.9 & \tabmaxcols & 100.0 &       & 2     & 105.7 & 40.7  & \tabmaxcols & 100.0 &       & 1     & 103.6 & 183.7 & 6,732,081 & 100.0 &       & 1     &       & 168.9 \\
A-n34-k5a & \tabopt & 39.7  & 60,015 & 100.0 &       & 1     &       & 1.5   &       &       &       &       &       &       &       &       &       &       &       &  \\
A-n34-k5b & \tabopt & 35.5  & 460   & 100.0 &       & 2     &       & 0.0   &       &       &       &       &       &       &       &       &       &       &       &  \\
A-n36-k5a & \tabopt & 151.8 & \tabmaxcols & 100.0 &       & 2     &       & 7.6   &       &       &       &       &       &       &       &       &       &       &       &  \\
A-n36-k5b & \tabopt & 165.7 & 74,972 & 100.0 &       & 2     &       & 3.2   &       &       &       &       &       &       &       &       &       &       &       &  \\
A-n37-k5a & \tabopt & 98.9  & 290,301 & 100.0 &       & 2     &       & 37.0  &       &       &       &       &       &       &       &       &       &       &       &  \\
A-n37-k5b & \tabopt & 166.0 & \tabmaxcols & 100.0 &       & 2     & 102.0 & 31.7  & 1,401,121 & 100.0 &       & 1     &       & 48.7  &       &       &       &       &       &  \\
A-n37-k6a & \tabopt & 27.0  & 171,702 & 100.0 &       & 2     &       & 3.5   &       &       &       &       &       &       &       &       &       &       &       &  \\
A-n37-k6b & \tabopt & 31.2  & \tabmaxcols & 100.0 &       & 2     & 100.2 & 5.5   & 259,368 & 100.0 &       & 1     &       & 2.2   &       &       &       &       &       &  \\
A-n38-k5a & \tabopt & 99.4  & \tabmaxcols & 100.0 &       & 2     & 101.1 & 5.7   & 352,959 & 100.0 &       & 1     &       & 32.4  &       &       &       &       &       &  \\
A-n38-k5b & \tabopt & 64.8  & \tabmaxcols & 100.0 &       & 2     & 103.0 & 4.9   & 1,152,405 & 100.0 &       & 1     &       & 8.2   &       &       &       &       &       &  \\
A-n39-k5a & \tabopt & 155.8 & \tabmaxcols & 100.0 &       & 2     & 102.4 & 7.2   & \tabmaxcols & 100.0 &       & 1     & 100.8 & 12.1  & 2,701,913 & 100.0 &       & 1     &       & 21.6 \\
A-n39-k5b & \tabopt & 182.7 & \tabmaxcols & 100.0 &       & 2     & 101.9 & 9.5   & \tabmaxcols & 100.0 &       & 1     & 100.3 & 15.0  & 1,909,058 & 100.0 &       & 1     &       & 19.1 \\
A-n39-k6a & \tabopt & 222.9 & \tabmaxcols & 99.7  &       & 2     & 103.3 & 8.8   & \tabmaxcols & 100.0 &       & 2     & 100.8 & 65.3  & 1,997,245 & 100.0 &       & 1     &       & 41.8 \\
A-n39-k6b & \tabopt & 100.6 & \tabmaxcols & 100.0 &       & 2     & 102.7 & 6.8   & \tabmaxcols & 100.0 &       & 1     & 100.7 & 12.9  & 2,070,281 & 100.0 &       & 1     &       & 16.7 \\
A-n44-k6a & \tabopt & 90.0  & \tabmaxcols & 100.0 &       & 2     & 102.4 & 8.6   & 1,088,748 & 100.0 &       & 1     &       & 11.7  &       &       &       &       &       &  \\
A-n44-k6b & \tabnoopt & 286.3 & \tabmaxcols & 100.0 &       & 2     & 102.9 & 29.5  & \tabmaxcols & 100.0 &       & 1     & 102.0 & 33.1  & \tabmaxcols & 100.0 &       & 1     & 100.7 & 126.9 \\
A-n45-k6a & \tabopt & 202.9 & \tabmaxcols & 100.0 &       & 2     &       & 132.9 &       &       &       &       &       &       &       &       &       &       &       &  \\
A-n45-k6b & \tabopt & 153.5 & \tabmaxcols & 100.0 &       & 2     & 101.3 & 18.1  & 1,220,510 & 100.0 &       & 1     &       & 70.7  &       &       &       &       &       &  \\
A-n45-k7a & \tabopt & 73.5  & \tabmaxcols & 100.0 &       & 2     & 105.3 & 6.9   & \tabmaxcols & 100.0 &       & 1     & 103.0 & 12.9  & 3,416,812 & 100.0 &       & 1     &       & 26.6 \\
A-n45-k7b & \tabopt & 75.2  & \tabmaxcols & 100.0 &       & 2     & 104.1 & 6.8   & \tabmaxcols & 100.0 &       & 1     & 102.4 & 12.4  & 3,568,626 & 100.0 &       & 1     &       & 27.1 \\
A-n46-k7a & \tabopt & 102.8 & \tabmaxcols & 100.0 &       & 2     & 101.8 & 39.9  & 970,793 & 100.0 &       & 1     &       & 12.0  &       &       &       &       &       &  \\
A-n46-k7b & \tabnoopt & 235.0 & \tabmaxcols & 100.0 &       & 2     & 104.0 & 15.8  & \tabmaxcols & 100.0 &       & 1     & 102.5 & 24.7  & \tabmaxcols & 100.0 &       & 1     & 100.5 & 131.2 \\
A-n48-k7a & \tabopt & 82.2  & 226,299 & 100.0 &       & 2     &       & 22.9  &       &       &       &       &       &       &       &       &       &       &       &  \\
A-n48-k7b & \tabopt & 87.0  & \tabmaxcols & 100.0 &       & 2     & 100.3 & 7.4   & 278,742 & 100.0 &       & 1     &       & 17.3  &       &       &       &       &       &  \\
A-n53-k7a & \tabopt & 503.7 & \tabmaxcols & 98.8  &       & 2     & 103.5 & 40.3  & \tabmaxcols & 100.0 &       & 2     & 101.1 & 225.1 & 2,185,209 & 100.0 &       & 1     &       & 166.9 \\
A-n53-k7b & \tabnoopt & 355.6 & \tabmaxcols & 99.8  &       & 2     & 104.1 & 21.8  & \tabmaxcols & 100.0 &       & 2     & 103.2 & 66.6  & \tabmaxcols & 100.0 &       & 1     & 102.1 & 183.1 \\
A-n54-k7a & \tabnoopt & 265.7 & \tabmaxcols & 100.0 &       & 2     & 102.7 & 8.0   & \tabmaxcols & 100.0 &       & 1     & 101.7 & 15.7  & \tabmaxcols & 100.0 &       & 1     & 100.4 & 59.4 \\
A-n54-k7b & \tabnoopt & 1065.3 & \tabmaxcols & 99.7  &       & 2     & 104.8 & 146.6 & \tabmaxcols & 100.0 &       & 2     & 104.1 & 414.6 & \tabmaxcols & 100.0 &       & 1     & 103.2 & 412.3 \\
A-n55-k9a & \tabnoopt & 599.1 & \tabmaxcols & 96.4  &       & 2     & 109.6 & 78.4  & \tabmaxcols & 97.7  &       & 2     & 108.8 & 41.0  & \tabmaxcols & 100.0 &       & 3     & 107.8 & 372.5 \\
A-n55-k9b & \tabnoopt & 377.7 & \tabmaxcols & 97.5  &       & 1     & 106.5 & 4.3   & \tabmaxcols & 100.0 &       & 2     & 105.7 & 57.3  & \tabmaxcols & 100.0 &       & 1     & 104.9 & 204.3 \\
A-n60-k9a & \tabnoopt & 948.9 & \tabmaxcols & 98.1  &       & 2     & 106.6 & 54.8  & \tabmaxcols & 99.8  &       & 2     & 105.7 & 336.2 & \tabmaxcols & 100.0 &       & 2     & 104.5 & 444.7 \\
A-n60-k9b & \tabnoopt & 21619.5 & \tabmaxcols & 96.9  &       & 2     & 104.9 & 100.8 & \tabmaxcols & 96.9  &       & 1     & 104.3 & 2269.1 & \tabmaxcols & 100.0 &       & 3     & 103.5 & 19142.3 \\
A-n61-k9a & \tabnoopt & 592.5 & \tabmaxcols & 99.3  &       & 1     & 108.0 & 3.9   & \tabmaxcols & 99.3  &       & 1     & 107.1 & 24.5  & \tabmaxcols & 100.0 &       & 2     & 106.1 & 288.9 \\
A-n61-k9b & \tabnoopt & 1020.9 & \tabmaxcols & 98.2  &       & 2     & 106.2 & 32.6  & \tabmaxcols & 98.3  &       & 2     & 105.6 & 366.7 & \tabmaxcols & 100.0 &       & 2     & 104.9 & 469.8 \\
A-n62-k8a & \tabnoopt & 446.2 & \tabmaxcols & 95.5  &       & 2     & 104.9 & 9.6   & \tabmaxcols & 100.0 &       & 2     & 104.0 & 54.5  & \tabmaxcols & 100.0 &       & 1     & 103.1 & 122.0 \\
A-n62-k8b & \tabnoopt & 1432.0 & \tabmaxcols & 99.5  &       & 2     & 104.8 & 21.4  & \tabmaxcols & 100.0 &       & 2     & 104.3 & 572.1 & \tabmaxcols & 100.0 &       & 1     & 103.7 & 562.1 \\
A-n63-k10a & \tabnoopt & 1446.4 & \tabmaxcols & 98.8  &       & 2     & 107.0 & 16.7  & \tabmaxcols & 99.9  &       & 2     & 106.0 & 1043.6 & \tabmaxcols & 100.0 &       & 2     & 105.1 & 313.5 \\
A-n63-k10b & \tabnoopt & 348.7 & \tabmaxcols & 99.1  &       & 2     & 104.7 & 61.8  & \tabmaxcols & 100.0 &       & 2     & 104.1 & 61.1  & \tabmaxcols & 100.0 &       & 1     & 103.4 & 132.7 \\
A-n63-k9a & \tabopt & 77.7  & \tabmaxcols & 100.0 &       & 2     & 101.6 & 6.9   & 700,337 & 100.0 &       & 1     &       & 6.6   &       &       &       &       &       &  \\
A-n63-k9b & \tabnoopt & 242.9 & \tabmaxcols & 100.0 &       & 2     & 103.2 & 12.4  & \tabmaxcols & 100.0 &       & 1     & 102.3 & 27.6  & \tabmaxcols & 100.0 &       & 1     & 100.6 & 106.0 \\
A-n64-k9a & \tabnoopt & 894.1 & \tabmaxcols & 98.6  &       & 2     & 103.9 & 73.0  & \tabmaxcols & 99.6  &       & 2     & 103.2 & 87.4  & \tabmaxcols & 100.0 &       & 2     & 102.1 & 624.6 \\
A-n64-k9b & \tabnoopt & 1120.9 & \tabmaxcols & 98.0  &       & 2     & 104.7 & 37.7  & \tabmaxcols & 99.7  &       & 2     & 104.2 & 229.1 & \tabmaxcols & 100.0 &       & 2     & 103.7 & 677.5 \\
A-n65-k9a & \tabnoopt & 340.5 & \tabmaxcols & 99.4  &       & 2     & 105.2 & 8.0   & \tabmaxcols & 99.4  &       & 1     & 104.3 & 20.5  & \tabmaxcols & 100.0 &       & 2     & 103.0 & 209.2 \\
A-n65-k9b & \tabnoopt & 298.4 & \tabmaxcols & 99.1  &       & 1     & 106.0 & 3.9   & \tabmaxcols & 100.0 &       & 2     & 105.0 & 67.8  & \tabmaxcols & 100.0 &       & 1     & 103.8 & 127.1 \\
A-n69-k9a & \tabnoopt & 3493.2 & \tabmaxcols & 98.8  &       & 1     & 107.5 & 1.3   & \tabmaxcols & 99.7  &       & 2     & 105.4 & 754.6 & \tabmaxcols & 100.0 &       & 2     & 103.9 & 2579.9 \\
A-n69-k9b & \tabnoopt & 1833.4 & \tabmaxcols & 96.8  &       & 2     & 107.2 & 11.2  & \tabmaxcols & 97.2  &       & 2     & 106.4 & 161.4 & \tabmaxcols & 100.0 &       & 2     & 105.7 & 1559.4 \\
A-n80-k10a & \tabnoopt & 599.5 & \tabmaxcols & 98.4  &       & 2     & 104.1 & 36.1  & \tabmaxcols & 100.0 &       & 2     & 103.6 & 100.7 & \tabmaxcols & 100.0 &       & 1     & 102.9 & 263.9 \\
A-n80-k10b & \tabnoopt & 1240.5 & \tabmaxcols & 99.3  &       & 2     & 102.8 & 30.2  & \tabmaxcols & 99.8  &       & 2     & 102.4 & 334.2 & \tabmaxcols & 100.0 &       & 2     & 102.0 & 671.4 \\
\hline
\end{longtable}%
}
}
\end{landscape}


\clearpage

\section{Conclusions}\label{sec:conclusions}

In this paper, we considered vehicle routing problems that can be modelled as a Set Partitioning (\fsp) problem with a linear fractional objective function. More precisely, we considered two objective functions: minimization of cost over load (also known as logistic ratio) and maximization of profit over time.

We investigated both continuous and  integer relaxations of the \fsp\ model. In particular, we proposed an alternative transformation to the transformation proposed by \citet{CC62} for linear fractional programming and a dual ascent heuristic used to compute both dual and primal bounds. The dual and primal bounds computed are embedded in an iterative exact procedure where at each iteration a reduced \fsp\ problem is solved by the extension of Dinkelbach's algorithm for fractional programming to integer programs.

We reported computational results showing that the proposed method solves to optimality instances involving up to 79 customers. The method can be easily adapted to deal with other routing constraints, simply by taking into account of such constraints in the route generation phase.










\clearpage

\section*{Appendix}

\subsection*{Proofs of statements}

\noindent {\bf Theorem 2.}
Let $\mb{v}$ be a feasible solution of problem $DNC\msp$ of cost $z(DNC\msp)$. A feasible solution $(\gb{\mu}, \omega)$ of $DCC\msp$ of cost $z($DCC\msp$)=z(DNC\msp)$ can be obtained by setting:
\begin{equation}\tag{\ref{thm2-a}}
  \omega=\sum_{i \in \setfc}v_i + \nveh v_0, \;\; \mu_0=\beta v_0, \;\; \mu_i=\beta v_i - s_i \omega, \forall i \in \setf, \;\; \mu_i=\beta v_i, \forall i \in \setc.
\end{equation}
\textit{Proof.}
It is easy to see that $z($DCC\msp$)=z(DNC\msp)$ due to the definition of $\omega$ in expressions \eqref{thm2-a}.
We have
\begin{equation}\label{thm2-b}
\begin{split}
  -\sum_{i \in \setfc} \mu_i -\nveh \mu_0  = -\sum_{i \in \setf}\beta v_i + \sum_{i \in \setf}s_i \omega - \sum_{i \in \setc}\beta v_i -  \nveh \beta v_0= \\
  -\beta (\sum_{i \in \setfc}v_i + \nveh v_0 ) + \beta (\sum_{i \in \setfc}v_i + \nveh v_0) = 0,
\end{split}
\end{equation}
showing that the dual constraint associated with variable $u$ of $CC\msp$ is satisfied.
Further,
\begin{equation}\label{thm2-c}
\begin{split}
\sum_{i \in \setfc}a_{i\ell}\mu_i + \mu_0 + w_\ell \omega = \sum_{i \in \setf}a_{il}\beta v_i - \sum_{i \in \setf}a_{il}s_i (\sum_{i \in \setfc}v_i + \nveh v_0) + \\ \sum_{i \in \setc}a_{i\ell}\beta v_i + \beta v_0 + w_\ell (\sum_{i \in \setfc}v_i + \nveh v_0) = \\ \sum_{i \in \setf}a_{il}\beta v_i + \sum_{i \in \setc}a_{i\ell}\beta v_i + \beta v_0 + (w_\ell-\sum_{i \in \setf}a_{il}s_i) (\sum_{i \in \setfc}v_i + \nveh v_0) = \\
\sum_{i \in \setfc}(\beta a_{il}+ \ov{w}_\ell)v_i  + (\beta + m \ov{w}_\ell)v_0 \leq c_\ell,
\end{split}
\end{equation}
showing that the dual constraint associated with variable $y_\ell$, $\ell \in \lc{R}$, of $CC\msp$ is satisfied.
$\,\Box$

\noindent {\bf Theorem 4.}
$\ov{r}=n(\ov{\mb{x}})/d(\ov{\mb{x}})=z(\msp(\lcb{R}))$ if, and only if, for the parametric problem $FP(r)$,  $z(FP(r))=\min \mml n(\mb{x})-r d(\mb{x}): \mb{x} \in P \mmr$,
$z(FP(\ov{r}))=n(\ov{\mb{x}})-\ov{r}d(\ov{\mb{x}})=0$.
\noindent {\textit{Proof.}}
\begin{enumerate}[(a)]
\item Let $\ov{\mb{x}}$ be an optimal solution of problem $\msp(\lcb{R})$. We have $
      \ov{r}=n(\ov{\mb{x}})/d(\ov{\mb{x}}) \leq n(\mb{x})/d(\mb{x}), \forall \mb{x} \in P$.
Hence
    \begin{enumerate}[(i)]
      \item $n(\mb{x}) - \ov{r}d(\mb{x}) \geq 0, \forall \mb{x} \in P$, and
      \item $n(\ov{\mb{x}})-\ov{r}d(\ov{\mb{x}})=0$.
    \end{enumerate}
This implies that $\ov{\mb{x}}$ is an optimal solution of $FP(\ov{r})$ of value $z(FP(\ov{r}))=0$.
\item Let $\mb{\ov{x}}$ be an optimal solution of $FP(\ov{r})$ such that $z(FP(\ov{r}))=n(\ov{\mb{x}})-\ov{r}d(\ov{\mb{x}})=0$. We have $n(\mb{x}) - \ov{r}d(\mb{x}) \geq 0, \forall \mb{x} \in P$, and $\ov{r} \leq n(\mb{x})/d(\mb{x}), \forall \mb{x} \in P$, therefore $\ov{r}$ is the minimum of problem $\msp(\lcb{R})$ that is taken at $\ov{\mb{x}}$.$\Box$
\end{enumerate}

\noindent {\bf Lemma 2.}
The following properties hold about the exact algorithm for solving $\msp(\lcb{R})$.
  \begin{enumerate}[(1)]
    \item The sequence $\{r_i\}$ is monotone decreasing, i.e., $r_i > r_{i+1}$, for all $i$ such that $z_i(FP(r_i)) <0$.
    \item \rev{If $z_{i}(FP(r_i))\geq0$ and solution $\mb{x}^i$ is optimal}, the value $\ov{z}(\msp(\lcb{R}))$ corresponds to the optimal solution value and $\ov{r}=n(\ov{\mb{x}})/d(\ov{\mb{x}})$.
    \item \rev{If $z_{i}(FP(r_i))<0$, we have} $d(\mb{x}^i) > d(\mb{x}^{i+1})$.
  \end{enumerate}
\noindent {\textit{Proof.}}
\begin{enumerate}[(1)]
    \item Since $d(\mb{x}^i)>0$ and $z_i(FP(r_i))=n(\mb{x}^i)-r_i d(\mb{x}^i)$ we have
        \begin{equation}\label{eq:lemmaexm-0}
          r_i > n(\mb{x}^i)/d(\mb{x}^i)=r_{i+1}.
        \end{equation}
    \item Assume that the algorithm terminates on the \rev{$i$-th iteration}. We have $z_i(FP(r_i)) \geq 0$. Consider a value $\hat{r}> \ov{r}=r_i$, we have:
        \begin{enumerate}[(i)]
          \item $\ov{r}=n(\mb{x}^{i-1})/d(\mb{x}^{i-1})$ and $n(\mb{x}^{i-1})-\ov{r}d(\mb{x}^{i-1})=0$
          \item As $d(\mb{x}^{i-1}) > 0$, we have                  $n(\mb{x}^{i-1})-\hat{r}d(\mb{x}^{i-1})<0$,
                  therefore value $\hat{r}$ does not satisfies the termination condition.
        \end{enumerate}
    \item We have:
        \begin{equation}\label{eq:lemmaexm}
            \left .
                \begin{array}{ll}
                (a)  \; &  n(\mb{x}^{i})-r_{i+1}d(\mb{x}^{i})=0 \; (\text{since the definition of } r_{i+1}),  \\
                (b)  \; &  n(\mb{x}^{i+1})-r_{i+1}d(\mb{x}^{i+1}) < 0 \; (z_i(FP(r_{i+1}))<0),  \\
                (c)  \; &  n(\mb{x}^{i+1})-r_i d(\mb{x}^{i+1}) \geq n(\mb{x}^{i})-r_i d(\mb{x}^{i}) \;(\mb{x}^{i+1} \text{ is a feasible solution of } FP(r_i)).
                \end{array}
            \right \}
        \end{equation}
        From \eqref{eq:lemmaexm}-a and \eqref{eq:lemmaexm}-b we have
        \begin{equation}\label{eq:lemmaexm-1}
          n(\mb{x}^{i+1}) - n(\mb{x}^{i}) \leq r_{i+1}(d(\mb{x}^{i+1})-d(\mb{x}^{i}))
        \end{equation}
        and by considering \eqref{eq:lemmaexm}-c from \eqref{eq:lemmaexm-1} we have
          $(r_{i+1}-r_i)(d(\mb{x}^{i+1})-d(\mb{x}^{i})) \geq 0$.
        From the above inequality and inequality \eqref{eq:lemmaexm-0} we have $d(\mb{x}^{i+1}) \leq d(\mb{x}^{i})$. Below we show that $d(\mb{x}^{i+1}) \neq d(\mb{x}^{i})$. Let $P_i=\{\mb{x} \in P: d(\mb{x})=d(\mb{x}^i)\}$ - we have $\mb{x}^i=\argmin \mml n(\mb{x})-r_i d(\mb{x}): \mb{x} \in P_i \mmr=\argmin \mml n(\mb{x}): \mb{x} \in P_i) \mmr$. If $\mb{x} \in P_i$, we have
          $n(\mb{x})-r_{i+1}d(\mb{x}) = n(\mb{x})-r_{i+1}d(\mb{x}^i) \geq n(\mb{x}^i)-r_{i+1}d(\mb{x}^i)=0$,
        therefore $n(\mb{x})-r_{i+1}d(\mb{x}) \geq 0$, $\forall \mb{x} \in P_i$.
        Since $r_{i+1}$ is not optimal, we have $z_{i+1}(FP(r_{i+1}))=n(\mb{x}^{i+1})-r_{i+1}d(\mb{x}^{i+1}) < 0$, therefore $\mb{x}^{i+1}$ is not in $P_i$. \;$\Box$
\end{enumerate}

\noindent {\bf Theorem 5.}
  For all $r_i \neq \ov{r}$ we have
  \begin{equation}\label{thm3-1}
    \frac{\ov{r}-r_{i+1}}{\ov{r}-r_i} \leq 1 - \frac{d(\ov{\mb{x}})}{d(\mb{x}^i)} < 1.
  \end{equation}
\noindent {\textit{Proof.}}
Inequality $n(\mb{x}^i)-r_i d(\mb{x}^i) \leq n(\ov{\mb{x}})-r_i d(\ov{\mb{x}})$ implies
\begin{equation}
  \frac{n(\mb{x}^i)}{d(\mb{x}^i)}-r_i \leq \frac{n(\ov{\mb{x}})}{d(\mb{x}^i)}-r_i\frac{d(\ov{\mb{x}})}{d(\mb{x}^i)},
\end{equation}
therefore
\begin{equation}\label{thm3-4}
    \begin{split}
    r_{i+1}-\ov{r} = \frac{n(\mb{x}^i)}{d(\mb{x}^i)} - \frac{n(\ov{\mb{x}})}{d(\ov{\mb{x}})} \leq
    \frac{n(\ov{\mb{x}})}{d(\mb{x}^i)}-\frac{n(\ov{\mb{x}})}{d(\ov{\mb{x}})}+r_i\left (1-\frac{d(\ov{\mb{x}})}{d(\mb{x}^i)} \right ) = \\
    \left ( \frac{1}{d(\mb{x}^i)}-\frac{1}{d(\ov{\mb{x}})} \right ) (n(\ov{\mb{x}})-r_id(\ov{\mb{x}}))) = \left ( \frac{1}{d(\mb{x}^i)}-\frac{1}{d(\ov{\mb{x}})} \right ) (\ov{r}-r_i)d(\ov{\mb{x}}),
    \end{split}
\end{equation}
since $\ov{r}=n(\ov{\mb{x}})/d(\ov{\mb{x}})$.
Since from Lemma \ref{lemmaexm}-(1) we have $\ov{r} < r_i$ and dividing by $(r_i-\ov{r}) > 0$ and since $d(\mb{x}^i)
 > d(\ov{\mb{x}})$ (Lemma \ref{lemmaexm}-(3)) we obtain
\begin{equation}\label{thm3-4}
    \frac{r_{i+1}-\ov{r}}{r_i-\ov{r}} \leq \left ( \frac{1}{d(\mb{x}^i)}-\frac{1}{d(\ov{\mb{x}})} \right ) \frac{\ov{r}-r_i}{r_i-\ov{r}}d(\ov{\mb{x}}) = \left ( - \frac{d(\ov{\mb{x}})}{d(\mb{x}^i)} + 1 \right ) < 1. \; \Box
\end{equation}

\subsection{Details about the computational results}

Tables \ref{tab:lb-CA}-\ref{tab:lb-PA} report the details about the bounding procedures. In particular,for bounding procedures $\cb$, $\ccg$ and \cg, the table reports the percentage deviation of the dual bound (column $\%B$) and the total computing time in seconds (column ``Time''). The percentage deviation is computed as $100.0 \er{\times} B/z^\ast$, where $z^\ast$ is the cost of the best solution found and $B$ is the value of the dual bound; for classes \classCA, $B$ refers to a \er{lower bound} whereas for class \classPA, $B$ refers to a \er{upper bound}. For bounding procedure \dk, the table reports only the computing time in seconds, being the value of the dual bound computed by the procedure equal to the one computed by procedure \cg.

For bounding procedures $\cb$ and $\ccg$, the table also reports the percentage deviation of the primal bound obtained by the heuristic procedure (column $\%PB$), computed as $100.0 \er{\times} PB/z^\ast$, where $z^\ast$ is the cost of the best solution found and $PB$ is the value of the primal bound - column ``Time'' also includes the time spent for computing the primal bound and, in the case of procedures $\cb$, the time spent for computing value $\tmin$.

For procedures \ccg, \cg, and \dk, the table also shows the number of columns or variables of the final master problem. In particular, for procedure \dk\ the number is computed as the sum over all algorithm iterations, whose average number is reported under column ``Iter'' in the table.

The percentage deviations are computed with respect to the values of the best solutions found which are reported in Tables \ref{tab:sol-CA}-\ref{tab:sol-PA}.

Tables \ref{tab:sol-CA}-\ref{tab:sol-PA} show details about the best solutions found. More precisely, for each instances, the tables show the values of the best solution found (\tabzbest) (including the corresponding numerator ``Cost'' or ``Profit'' and denominator ``Load'' or ``Time''), the number of routes of the solution (\tabsolm ), the number of optional customers selected in solution (\tabnfixed), the percentage of the working time utilization (\tabut) and the average number of customers per route (\tabavgc).



\begin{table}
\caption{Dual and primal bounds on instances of class \classCA} \label{tab:lb-CA}
\begin{center}
{\scriptsize{
\setlength{\tabcolsep}{2.7pt}
\renewcommand{\arraystretch}{1.0}
\begin{tabular}{r|rrr|rrrr|rrr|rrrr}
\hline
      & \multicolumn{3}{c|}{Procedure \cb} & \multicolumn{4}{c|}{Procedure \ccg}     & \multicolumn{3}{c|}{Procedure \cg} & \multicolumn{3}{c}{Procedure \dk} \bigstrut[t]\\
\tabname & \multicolumn{1}{c}{$\%B$} & \multicolumn{1}{c}{$\%PB$} & \multicolumn{1}{c|}{\tabtime} & \multicolumn{1}{c}{$\%B$} & \multicolumn{1}{c}{$\%PB$} & \multicolumn{1}{c}{\tablbcols} & \multicolumn{1}{c|}{\tabtime} & \multicolumn{1}{c}{$\%B$} & \multicolumn{1}{c}{\tablbcols} & \multicolumn{1}{c|}{\tabtime} & \multicolumn{1}{c}{\tablbit} & \multicolumn{1}{c}{\tablbcols} & \multicolumn{1}{c}{\tabtime} \bigstrut[b]\\
\hline
A-n32-k5a & 98.6  & 103.1 & 7.2   & 99.3  & 103.1 & 160   & 0.1   & 99.3  & 2,690 & 0.2   & 3     & 2,394 & 0.2 \bigstrut[t]\\
A-n32-k5b & 98.1  & 102.2 & 18.8  & 98.2  & 102.2 & 232   & 0.1   & 98.2  & 3,340 & 0.2   & 3     & 3,439 & 0.3 \\
A-n33-k5a & 98.1  & 112.5 & 18.4  & 98.1  & 112.5 & 252   & 0.1   & 98.1  & 2,689 & 0.1   & 2     & 2,503 & 0.2 \\
A-n33-k5b & 99.4  & 100.4 & 22.3  & 99.5  & 100.4 & 150   & 0.1   & 99.5  & 2,765 & 0.1   & 3     & 3,090 & 0.2 \\
A-n33-k6a & 99.5  & 119.4 & 6.3   & 100.0 & 119.4 & 147   & 0.1   & 100.0 & 2,079 & 0.1   & 3     & 1,787 & 0.1 \\
A-n33-k6b & 97.9  & 103.0 & 13.3  & 97.7  & 103.0 & 261   & 0.1   & 97.9  & 2,444 & 0.1   & 3     & 2,379 & 0.1 \\
A-n34-k5a & 99.4  & 107.5 & 19.6  & 99.4  & 107.5 & 207   & 0.1   & 99.4  & 2,731 & 0.1   & 3     & 3,237 & 0.2 \\
A-n34-k5b & 98.8  & 106.6 & 16.5  & 98.9  & 106.6 & 314   & 0.1   & 98.9  & 3,222 & 0.2   & 3     & 2,777 & 0.2 \\
A-n36-k5a & 99.5  & 115.4 & 25.2  & 99.7  & 113.4 & 315   & 0.5   & 99.7  & 4,103 & 0.4   & 3     & 4,356 & 0.8 \\
A-n36-k5b & 99.0  & 106.6 & 25.3  & 98.8  & 106.5 & 411   & 0.4   & 99.1  & 5,104 & 0.6   & 3     & 5,691 & 0.9 \\
A-n37-k5a & 96.1  & 102.8 & 17.0  & 97.5  & 102.8 & 200   & 0.2   & 97.5  & 4,020 & 0.4   & 3     & 4,126 & 0.4 \\
A-n37-k5b & 98.4  & 104.8 & 23.5  & 98.3  & 103.8 & 396   & 0.2   & 98.6  & 7,444 & 1.2   & 3     & 6,008 & 0.8 \\
A-n37-k6a & 99.2  & 105.3 & 8.1   & 99.3  & 103.6 & 313   & 0.1   & 99.3  & 2,899 & 0.1   & 3     & 2,595 & 0.2 \\
A-n37-k6b & 97.4  & 110.6 & 11.1  & 97.5  & 112.5 & 396   & 0.2   & 97.5  & 3,795 & 0.3   & 3     & 3,763 & 0.4 \\
A-n38-k5a & 99.8  & 112.5 & 4.7   & 100.0 & 112.5 & 299   & 0.1   & 100.0 & 3,915 & 0.3   & 2     & 4,434 & 0.3 \\
A-n38-k5b & 95.7  & 108.1 & 10.1  & 95.9  & 107.4 & 328   & 0.2   & 95.9  & 4,434 & 0.3   & 2     & 4,204 & 0.3 \\
A-n39-k5a & 97.8  & 106.4 & 33.5  & 97.8  & 109.2 & 313   & 0.2   & 97.8  & 3,646 & 0.3   & 3     & 3,799 & 0.5 \\
A-n39-k5b & 99.2  & 100.8 & 36.4  & 99.3  & 102.9 & 329   & 0.3   & 99.3  & 5,449 & 0.5   & 2     & 4,850 & 0.5 \\
A-n39-k6a & 97.8  & 107.9 & 11.0  & 98.2  & 107.9 & 197   & 0.2   & 98.2  & 3,505 & 0.2   & 3     & 3,133 & 0.2 \\
A-n39-k6b & 98.1  & 100.9 & 16.7  & 98.4  & 100.9 & 278   & 0.3   & 98.4  & 4,878 & 0.4   & 3     & 3,788 & 0.3 \\
A-n44-k6a & 100.0 & 118.9 & 10.0  & 100.0 & 118.9 & 349   & 0.2   & 100.0 & 4,661 & 0.4   & 2     & 5,211 & 0.5 \\
A-n44-k6b & 97.7  & 102.8 & 15.9  & 97.8  & 102.7 & 484   & 0.3   & 97.8  & 5,620 & 0.6   & 3     & 5,744 & 0.7 \\
A-n45-k6a & 99.0  & 101.7 & 13.7  & 99.1  & 102.5 & 220   & 0.4   & 99.1  & 4,738 & 0.4   & 3     & 4,647 & 0.5 \\
A-n45-k6b & 98.5  & 104.6 & 21.0  & 98.6  & 104.1 & 354   & 0.2   & 98.6  & 6,183 & 0.6   & 3     & 5,769 & 0.7 \\
A-n45-k7a & 98.7  & 104.0 & 12.2  & 98.8  & 102.7 & 304   & 0.2   & 98.8  & 4,597 & 0.4   & 3     & 5,268 & 0.5 \\
A-n45-k7b & 99.4  & 106.4 & 16.1  & 99.6  & 105.9 & 350   & 0.3   & 99.6  & 5,647 & 0.5   & 3     & 5,043 & 0.5 \\
A-n46-k7a & 98.3  & 103.5 & 9.6   & 98.7  & 100.7 & 221   & 0.3   & 98.7  & 5,387 & 0.5   & 3     & 4,206 & 0.4 \\
A-n46-k7b & 96.8  & 102.4 & 6.8   & 98.5  & 102.4 & 352   & 0.4   & 98.5  & 6,694 & 0.7   & 3     & 6,927 & 0.8 \\
A-n48-k7a & 97.6  & 106.1 & 12.7  & 97.7  & 106.1 & 214   & 0.2   & 97.7  & 4,831 & 0.4   & 3     & 5,146 & 0.5 \\
A-n48-k7b & 98.5  & 104.5 & 20.2  & 98.7  & 104.5 & 413   & 0.3   & 98.7  & 7,221 & 1.1   & 2     & 7,214 & 1.1 \\
A-n53-k7a & 100.0 & 116.0 & 21.7  & 100.0 & 105.7 & 405   & 0.2   & 100.0 & 9,810 & 1.4   & 2     & 9,232 & 1.4 \\
A-n53-k7b & 99.8  & 100.8 & 26.9  & 100.0 & 100.2 & 261   & 0.3   & 100.0 & 13,143 & 2.7  & 3     & 13,137 & 2.9 \\
A-n54-k7a & 99.7  & 111.1 & 15.1  & 99.9  & 111.1 & 426   & 0.3   & 99.9  & 7,265 & 0.9   & 3     & 6,080 & 0.8 \\
A-n54-k7b & 98.1  & 105.2 & 18.3  & 98.0  & 104.4 & 574   & 0.5   & 98.2  & 9,760 & 1.7   & 3     & 11,170 & 2.3 \\
A-n55-k9a & 98.5  & 108.5 & 11.6  & 98.8  & 107.9 & 379   & 0.4   & 98.8  & 4,514 & 0.3   & 3     & 4,487 & 0.4 \\
A-n55-k9b & 98.8  & 105.2 & 15.0  & 98.9  & 105.1 & 456   & 0.3   & 98.9  & 7,061 & 0.7   & 3     & 6,687 & 0.7 \\
A-n60-k9a & 98.2  & 110.2 & 23.6  & 98.4  & 102.7 & 264   & 0.6   & 98.4  & 9,161 & 1.2   & 2     & 8,787 & 0.9 \\
A-n60-k9b & 98.5  & 109.0 & 35.0  & 98.6  & 106.1 & 677   & 1.0   & 98.6  & 10,938 & 1.9  & 3     & 11,252 & 2.3 \\
A-n61-k9a & 98.0  & 109.8 & 14.9  & 98.2  & 109.9 & 468   & 0.4   & 98.2  & 6,989 & 0.7   & 3     & 6,510 & 0.7 \\
A-n61-k9b & 99.5  & 106.0 & 19.0  & 99.8  & 106.0 & 441   & 0.5   & 99.8  & 9,001 & 1.3   & 3     & 8,235 & 1.2 \\
A-n62-k8a & 97.7  & 103.3 & 13.6  & 98.6  & 100.5 & 604   & 1.0   & 98.6  & 13,514 & 2.2  & 3     & 12,301 & 2.1 \\
A-n62-k8b & 98.0  & 102.9 & 27.1  & 98.2  & 102.9 & 502   & 1.0   & 98.2  & 15,392 & 3.4  & 3     & 18,403 & 4.4 \\
A-n63-k10a & 98.8  & 111.3 & 18.0  & 98.9  & 110.3 & 296   & 0.5   & 98.9  & 8,630 & 1.1  & 2     & 7,417 & 0.8 \\
A-n63-k10b & 98.1  & 102.6 & 23.0  & 98.7  & 101.6 & 409   & 1.1   & 98.7  & 11,353 & 1.7 & 3     & 10,505 & 1.6 \\
A-n63-k9a & 98.6  & 109.7 & 19.0  & 98.7  & 108.1 & 369   & 0.4   & 98.7  & 11,333 & 1.6  & 2     & 9,417 & 1.3 \\
A-n63-k9b & 98.6  & 104.2 & 30.2  & 98.8  & 104.2 & 494   & 0.7   & 98.8  & 11,169 & 1.7  & 3     & 9,443 & 1.7 \\
A-n64-k9a & 99.1  & 108.6 & 12.9  & 99.3  & 109.3 & 368   & 0.6   & 99.3  & 10,926 & 1.5  & 3     & 10,341 & 1.4 \\
A-n64-k9b & 97.5  & 112.0 & 8.6   & 99.1  & 106.8 & 719   & 1.0   & 99.1  & 15,511 & 3.2  & 2     & 14,823 & 3.4 \\
A-n65-k9a & 98.4  & 104.3 & 17.0  & 98.8  & 111.6 & 344   & 0.3   & 98.8  & 8,175 & 1.0   & 3     & 7,978 & 1.0 \\
A-n65-k9b & 98.5  & 104.1 & 27.4  & 98.6  & 106.2 & 404   & 0.7   & 98.6  & 14,128 & 2.5  & 3     & 11,708 & 2.2 \\
A-n69-k9a & 98.9  & 112.4 & 34.1  & 99.1  & 110.9 & 370   & 0.8   & 99.1  & 15,303 & 3.0  & 3     & 13,363 & 2.4 \\
A-n69-k9b & 97.5  & 110.7 & 49.3  & 97.8  & 107.8 & 675   & 1.0   & 97.8  & 17,575 & 4.3  & 3     & 15,376 & 3.4 \\
A-n80-k10a & 99.2  & 108.7 & 18.3  & 99.4  & 107.9 & 775   & 1.3   & 99.4  & 16,459 & 2.8 & 3     & 15,720 & 3.0 \\
A-n80-k10b & 98.4  & 107.6 & 30.7  & 98.5  & 106.5 & 769   & 1.9   & 98.5  & 27,078 & 8.5 & 3     & 21,825 & 5.5 \bigstrut[b]\\
\hline
      & 98.5  & 107.0 & 18.8  & 98.7  & 106.3 & 374   & 0.4   & 98.8  & 7,684 & 1.2   & 2.8   & 7,254 & 1.1 \bigstrut\\
\hline
\end{tabular}%

}}
\end{center}
\end{table}

\begin{table}
\caption{Dual and primal bounds on instances of class \classPA} \label{tab:lb-PA}
\begin{center}
{\scriptsize{
\setlength{\tabcolsep}{2.7pt}
\renewcommand{\arraystretch}{1.0}
\begin{tabular}{r|rrr|rrrr|rrr|rrrr}
\hline
      & \multicolumn{3}{c|}{Procedure \cb} & \multicolumn{4}{c|}{Procedure \ccg}     & \multicolumn{3}{c|}{Procedure \cg} & \multicolumn{3}{c}{Procedure \dk} \bigstrut[t]\\
\tabname & \multicolumn{1}{c}{$\%B$} & \multicolumn{1}{c}{$\%PB$} & \multicolumn{1}{c|}{\tabtime} & \multicolumn{1}{c}{$\%B$} & \multicolumn{1}{c}{$\%PB$} & \multicolumn{1}{c}{\tablbcols} & \multicolumn{1}{c|}{\tabtime} & \multicolumn{1}{c}{$\%B$} & \multicolumn{1}{c}{\tablbcols} & \multicolumn{1}{c|}{\tabtime} & \multicolumn{1}{c}{\tablbit} & \multicolumn{1}{c}{\tablbcols} & \multicolumn{1}{c}{\tabtime} \bigstrut[b]\\
\hline
A-n32-k5a & 105.8 & 99.8  & 134.8 & 103.7 & 99.8  & 228   & 0.1   & 103.5 & 719   & 0.6   & 2     & 597   & 1.0 \bigstrut[t]\\
A-n32-k5b & 106.5 & 96.3  & 49.8  & 106.2 & 98.3  & 286   & 0.1   & 106.0 & 1,009 & 0.3   & 2     & 950   & 0.2 \\
A-n33-k5a & 100.5 & 99.3  & 69.0  & 100.5 & 99.9  & 413   & 0.2   & 100.3 & 1,367 & 0.4   & 2     & 1,874 & 1.9 \\
A-n33-k5b & 101.5 & 99.6  & 72.5  & 101.2 & 100.0 & 511   & 0.3   & 101.1 & 1,589 & 0.4   & 3     & 1,631 & 0.3 \\
A-n33-k6a & 110.7 & 97.7  & 47.9  & 110.6 & 98.3  & 586   & 0.1   & 110.4 & 1,412 & 0.2   & 3     & 1,754 & 1.0 \\
A-n33-k6b & 109.6 & 99.1  & 50.1  & 109.6 & 99.6  & 566   & 0.2   & 109.4 & 1,713 & 0.2   & 3     & 1,609 & 0.1 \\
A-n34-k5a & 107.5 & 98.2  & 37.5  & 106.8 & 100.0 & 323   & 0.1   & 106.8 & 1,366 & 0.3   & 3     & 1,292 & 0.6 \\
A-n34-k5b & 100.9 & 99.7  & 34.3  & 100.8 & 99.7  & 317   & 0.1   & 100.6 & 1,128 & 0.3   & 3     & 1,457 & 0.3 \\
A-n36-k5a & 105.2 & 94.2  & 141.2 & 104.7 & 95.7  & 458   & 0.4   & 104.6 & 1,133 & 0.4   & 2     & 1,051 & 0.8 \\
A-n36-k5b & 103.7 & 96.0  & 159.6 & 103.5 & 97.2  & 516   & 0.5   & 103.5 & 1,469 & 1.1   & 2     & 1,127 & 0.6 \\
A-n37-k5a & 111.1 & 98.8  & 61.2  & 110.2 & 98.0  & 464   & 0.1   & 110.0 & 1,292 & 0.2   & 2     & 1,150 & 0.4 \\
A-n37-k5b & 105.3 & 99.5  & 81.6  & 105.1 & 99.6  & 441   & 0.2   & 105.1 & 1,684 & 0.2   & 3     & 1,517 & 0.2 \\
A-n37-k6a & 104.1 & 97.8  & 22.4  & 103.6 & 97.8  & 522   & 0.3   & 103.4 & 1,281 & 0.1   & 2     & 1,377 & 0.2 \\
A-n37-k6b & 105.2 & 96.0  & 20.7  & 104.2 & 96.8  & 540   & 0.3   & 104.0 & 1,210 & 0.2   & 2     & 1,333 & 0.1 \\
A-n38-k5a & 106.5 & 95.3  & 61.5  & 106.3 & 96.4  & 395   & 0.1   & 106.2 & 1,406 & 0.2   & 3     & 1,641 & 0.2 \\
A-n38-k5b & 107.6 & 96.4  & 49.2  & 107.7 & 99.0  & 394   & 0.2   & 107.4 & 1,545 & 0.2   & 3     & 1,597 & 0.1 \\
A-n39-k5a & 105.9 & 98.0  & 106.7 & 106.0 & 96.1  & 371   & 0.3   & 105.8 & 1,349 & 0.5   & 3     & 1,970 & 1.0 \\
A-n39-k5b & 105.2 & 95.8  & 131.7 & 104.8 & 97.7  & 342   & 0.4   & 104.7 & 1,458 & 0.5   & 3     & 1,628 & 0.3 \\
A-n39-k6a & 107.6 & 98.6  & 101.6 & 107.4 & 98.6  & 450   & 0.2   & 107.4 & 1,140 & 0.6   & 3     & 1,584 & 0.6 \\
A-n39-k6b & 106.3 & 94.5  & 58.1  & 106.3 & 98.3  & 402   & 0.3   & 106.0 & 1,282 & 0.2   & 3     & 1,211 & 0.3 \\
A-n44-k6a & 109.5 & 94.4  & 65.8  & 107.6 & 97.0  & 545   & 0.4   & 107.5 & 1,568 & 0.6   & 2     & 1,485 & 0.5 \\
A-n44-k6b & 106.4 & 95.0  & 79.9  & 105.8 & 96.7  & 536   & 0.5   & 105.7 & 1,999 & 0.8   & 3     & 1,881 & 0.3 \\
A-n45-k6a & 102.5 & 95.9  & 67.1  & 102.5 & 97.3  & 611   & 0.3   & 102.3 & 1,980 & 0.5   & 3     & 2,056 & 0.3 \\
A-n45-k6b & 104.0 & 96.2  & 60.9  & 103.9 & 99.0  & 634   & 0.5   & 103.7 & 2,127 & 0.4   & 3     & 2,272 & 0.4 \\
A-n45-k7a & 110.9 & 97.7  & 19.7  & 110.0 & 98.8  & 461   & 0.3   & 109.8 & 1,168 & 0.1   & 2     & 1,142 & 0.1 \\
A-n45-k7b & 110.5 & 96.7  & 19.9  & 107.0 & 97.6  & 436   & 0.4   & 107.0 & 1,366 & 0.1   & 3     & 1,406 & 0.1 \\
A-n46-k7a & 106.5 & 98.9  & 47.9  & 106.1 & 99.2  & 547   & 0.3   & 105.8 & 2,146 & 0.3   & 2     & 1,849 & 0.7 \\
A-n46-k7b & 107.4 & 96.6  & 46.7  & 107.0 & 97.4  & 691   & 0.4   & 106.8 & 2,198 & 0.3   & 3     & 2,471 & 0.2 \\
A-n48-k7a & 103.4 & 98.6  & 56.7  & 103.6 & 99.0  & 411   & 0.3   & 103.3 & 1,473 & 0.4   & 3     & 1,641 & 0.3 \\
A-n48-k7b & 103.7 & 97.4  & 58.3  & 103.5 & 97.8  & 510   & 0.6   & 103.4 & 1,984 & 0.3   & 2     & 1,931 & 0.3 \\
A-n53-k7a & 107.7 & 90.4  & 63.1  & 107.0 & 96.5  & 523   & 0.3   & 107.0 & 2,529 & 0.5   & 2     & 2,049 & 1.1 \\
A-n53-k7b & 106.4 & 97.3  & 65.5  & 106.0 & 98.6  & 604   & 0.7   & 105.8 & 2,504 & 0.4   & 3     & 2,688 & 0.2 \\
A-n54-k7a & 105.7 & 95.2  & 163.6 & 105.4 & 97.2  & 596   & 0.8   & 105.2 & 2,139 & 1.6   & 3     & 2,281 & 2.2 \\
A-n54-k7b & 107.2 & 98.3  & 72.2  & 106.3 & 96.7  & 597   & 0.8   & 106.2 & 2,463 & 0.5   & 3     & 2,375 & 0.4 \\
A-n55-k9a & 111.5 & 91.5  & 91.6  & 110.7 & 95.3  & 513   & 0.4   & 110.5 & 2,421 & 0.7   & 3     & 3,156 & 0.6 \\
A-n55-k9b & 108.5 & 92.9  & 93.3  & 107.3 & 97.5  & 799   & 0.7   & 107.1 & 3,726 & 1.9   & 3     & 3,595 & 0.6 \\
A-n60-k9a & 108.9 & 96.1  & 95.7  & 108.5 & 98.6  & 681   & 1.2   & 108.4 & 2,892 & 1.1   & 3     & 3,157 & 2.2 \\
A-n60-k9b & 106.8 & 93.5  & 85.1  & 105.9 & 96.3  & 703   & 1.3   & 105.7 & 3,375 & 1.3   & 2     & 3,680 & 2.2 \\
A-n61-k9a & 110.9 & 98.2  & 258.5 & 109.4 & 99.8  & 751   & 1.1   & 109.2 & 2,854 & 3.4   & 3     & 4,223 & 2.0 \\
A-n61-k9b & 108.7 & 95.9  & 131.9 & 107.4 & 95.5  & 993   & 1.2   & 107.3 & 4,317 & 1.4   & 3     & 4,053 & 0.9 \\
A-n62-k8a & 106.7 & 91.1  & 240.2 & 106.0 & 94.2  & 813   & 1.4   & 105.9 & 3,140 & 1.2   & 2     & 3,101 & 1.9 \\
A-n62-k8b & 106.0 & 94.4  & 242.3 & 105.6 & 95.1  & 847   & 2.0   & 105.6 & 4,475 & 1.4   & 3     & 4,020 & 0.7 \\
A-n63-k10a & 109.2 & 96.1  & 51.9  & 108.7 & 98.4  & 764   & 0.6   & 108.6 & 3,072 & 0.8  & 3     & 3,730 & 1.2 \\
A-n63-k10b & 107.0 & 98.6  & 68.7  & 106.2 & 98.7  & 1,026 & 1.3   & 106.1 & 3,054 & 0.7  & 3     & 3,686 & 0.4 \\
A-n63-k9a & 105.8 & 98.0  & 56.4  & 105.8 & -     & 560   & 0.8   & 105.6 & 1,896 & 0.4   & 2     & 1,719 & 0.3 \\
A-n63-k9b & 109.7 & 96.5  & 72.7  & 104.9 & -     & 650   & 1.2   & 104.8 & 2,110 & 0.4   & 3     & 2,132 & 0.3 \\
A-n64-k9a & 107.6 & 95.1  & 88.9  & 105.6 & 95.6  & 812   & 1.5   & 105.5 & 2,176 & 0.7   & 3     & 2,423 & 0.8 \\
A-n64-k9b & 106.5 & 93.0  & 152.6 & 105.8 & -     & 1,325 & 1.5   & 105.7 & 4,154 & 1.6   & 3     & 4,439 & 2.2 \\
A-n65-k9a & 107.5 & 93.9  & 82.7  & 106.8 & 97.4  & 624   & 0.7   & 106.8 & 2,861 & 0.6   & 3     & 3,671 & 0.5 \\
A-n65-k9b & 107.7 & 93.6  & 76.9  & 106.9 & 96.0  & 746   & 1.2   & 106.8 & 3,185 & 0.6   & 3     & 3,283 & 0.4 \\
A-n69-k9a & 108.5 & 95.9  & 104.8 & 107.7 & 98.8  & 814   & 0.9   & 107.5 & 4,513 & 2.1   & 2     & 3,863 & 2.5 \\
A-n69-k9b & 109.1 & 94.5  & 78.1  & 108.0 & 96.0  & 951   & 1.3   & 107.9 & 5,028 & 1.4   & 3     & 4,633 & 0.7 \\
A-n80-k10a & 107.8 & 97.5  & 173.5 & 105.3 & 98.0  & 1,096 & 1.5   & 105.2 & 4,087 & 1.8  & 2     & 4,408 & 2.3 \\
A-n80-k10b & 104.9 & 95.9  & 175.2 & 103.6 & 97.4  & 1,250 & 2.6   & 103.6 & 5,697 & 2.6  & 2     & 5,535 & 2.4 \bigstrut[b]\\
\hline
      & 106.8 & 96.3  & 88.9  & 106.1 & 97.7  & 610   & 0.7   & 105.9 & 2,282 & 0.7   & 2.6   & 2,377 & 0.8 \bigstrut\\
\hline
\end{tabular}%

}}
\end{center}
\end{table}



\begin{table}
\caption{Details of the best solutions found for class \classCA} \label{tab:sol-CA}
\begin{center}
{\scriptsize{
\setlength{\tabcolsep}{2.7pt}
\renewcommand{\arraystretch}{1.0}
\begin{tabular}{rrrrr|rrrrrrr}
\hline
\tabname & $\setv$ & $\nf$ & $\nc$ & $\nveh$ & \tabzbest & Cost  & Load  & \tabsolm & \tabnfixed & \tabut & \tabavgc \bigstrut\\
\hline
A-n32-k5a & 32    & 15    & 16    & 4     & 1.8264 & 705   & 386   & 4     & 13    & 96.5  & 7.0 \bigstrut[t]\\
A-n32-k5b & 32    & 23    & 8     & 4     & 1.8677 & 734   & 393   & 4     & 6     & 98.3  & 7.3 \\
A-n33-k5a & 33    & 16    & 16    & 4     & 1.3984 & 523   & 374   & 4     & 11    & 93.5  & 6.8 \\
A-n33-k5b & 33    & 24    & 8     & 5     & 1.3960 & 557   & 399   & 5     & 3     & 79.8  & 5.4 \\
A-n33-k6a & 33    & 16    & 16    & 4     & 1.1654 & 444   & 381   & 4     & 5     & 95.3  & 5.3 \\
A-n33-k6b & 33    & 24    & 8     & 6     & 1.2990 & 630   & 485   & 6     & 4     & 80.8  & 4.7 \\
A-n34-k5a & 34    & 16    & 17    & 4     & 1.5106 & 500   & 331   & 4     & 8     & 82.8  & 6.0 \\
A-n34-k5b & 34    & 24    & 9     & 5     & 1.5606 & 682   & 437   & 5     & 5     & 87.4  & 5.8 \\
A-n36-k5a & 36    & 17    & 18    & 4     & 1.7247 & 664   & 385   & 4     & 12    & 96.3  & 7.3 \\
A-n36-k5b & 36    & 26    & 9     & 5     & 1.7532 & 689   & 393   & 4     & 5     & 98.3  & 7.8 \\
A-n37-k5a & 37    & 18    & 18    & 3     & 1.7508 & 520   & 297   & 3     & 7     & 99.0  & 8.3 \\
A-n37-k5b & 37    & 27    & 9     & 4     & 1.6292 & 637   & 391   & 4     & 6     & 97.8  & 8.3 \\
A-n37-k6a & 37    & 18    & 18    & 4     & 1.6051 & 634   & 395   & 4     & 6     & 98.8  & 6.0 \\
A-n37-k6b & 37    & 27    & 9     & 5     & 1.6862 & 833   & 494   & 5     & 4     & 98.8  & 6.2 \\
A-n38-k5a & 38    & 18    & 19    & 4     & 1.3813 & 547   & 396   & 4     & 10    & 99.0  & 7.0 \\
A-n38-k5b & 38    & 27    & 10    & 5     & 1.4846 & 677   & 456   & 5     & 7     & 91.2  & 6.8 \\
A-n39-k5a & 39    & 19    & 19    & 4     & 1.4987 & 559   & 373   & 4     & 10    & 93.3  & 7.3 \\
A-n39-k5b & 39    & 28    & 10    & 5     & 1.5894 & 658   & 414   & 5     & 4     & 82.8  & 6.4 \\
A-n39-k6a & 39    & 19    & 19    & 5     & 1.4158 & 664   & 469   & 5     & 12    & 93.8  & 6.2 \\
A-n39-k6b & 39    & 28    & 10    & 6     & 1.4969 & 723   & 483   & 6     & 4     & 80.5  & 5.3 \\
A-n44-k6a & 44    & 21    & 22    & 4     & 1.6286 & 627   & 385   & 4     & 9     & 96.3  & 7.5 \\
A-n44-k6b & 44    & 32    & 11    & 6     & 1.6313 & 907   & 556   & 6     & 10    & 92.7  & 7.0 \\
A-n45-k6a & 45    & 22    & 22    & 4     & 1.6500 & 627   & 380   & 4     & 6     & 95.0  & 7.0 \\
A-n45-k6b & 45    & 33    & 11    & 6     & 1.5650 & 867   & 554   & 6     & 7     & 92.3  & 6.7 \\
A-n45-k7a & 45    & 22    & 22    & 4     & 1.7051 & 665   & 390   & 4     & 7     & 97.5  & 7.3 \\
A-n45-k7b & 45    & 33    & 11    & 6     & 1.7082 & 960   & 562   & 6     & 6     & 93.7  & 6.5 \\
A-n46-k7a & 46    & 22    & 23    & 4     & 1.5291 & 604   & 395   & 4     & 7     & 98.8  & 7.3 \\
A-n46-k7b & 46    & 33    & 12    & 6     & 1.4663 & 827   & 564   & 6     & 8     & 94.0  & 6.8 \\
A-n48-k7a & 48    & 23    & 24    & 5     & 1.6511 & 814   & 493   & 5     & 12    & 98.6  & 7.0 \\
A-n48-k7b & 48    & 35    & 12    & 6     & 1.6712 & 986   & 590   & 6     & 9     & 98.3  & 7.3 \\
A-n53-k7a & 53    & 26    & 26    & 5     & 1.4429 & 720   & 499   & 5     & 10    & 99.8  & 7.2 \\
A-n53-k7b & 53    & 39    & 13    & 6     & 1.4661 & 821   & 560   & 6     & 5     & 93.3  & 7.3 \\
A-n54-k7a & 54    & 26    & 27    & 5     & 1.6430 & 810   & 493   & 5     & 12    & 98.6  & 7.6 \\
A-n54-k7b & 54    & 39    & 14    & 6     & 1.6965 & 1,006 & 593   & 6     & 7     & 98.8  & 7.7 \\
A-n55-k9a & 55    & 27    & 27    & 6     & 1.2271 & 724   & 590   & 6     & 12    & 98.3  & 6.5 \\
A-n55-k9b & 55    & 40    & 14    & 8     & 1.2105 & 897   & 741   & 8     & 6     & 92.6  & 5.8 \\
A-n60-k9a & 60    & 29    & 30    & 5     & 1.4888 & 734   & 493   & 5     & 7     & 98.6  & 7.2 \\
A-n60-k9b & 60    & 44    & 15    & 7     & 1.5768 & 1,099 & 697   & 7     & 8     & 99.6  & 7.4 \\
A-n61-k9a & 61    & 30    & 30    & 6     & 1.1142 & 634   & 569   & 6     & 10    & 94.8  & 6.7 \\
A-n61-k9b & 61    & 45    & 15    & 8     & 1.1207 & 882   & 787   & 8     & 7     & 98.4  & 6.5 \\
A-n62-k8a & 62    & 30    & 31    & 4     & 1.8693 & 744   & 398   & 4     & 6     & 99.5  & 9.0 \\
A-n62-k8b & 62    & 45    & 16    & 7     & 1.7125 & 1,120 & 654   & 7     & 8     & 93.4  & 7.6 \\
A-n63-k10a & 63    & 31    & 31    & 6     & 1.3548 & 798   & 589   & 6     & 10    & 98.2  & 6.8 \\
A-n63-k10b & 63    & 46    & 16    & 8     & 1.4383 & 1,142 & 794   & 8     & 8     & 99.3  & 6.8 \\
A-n63-k9a & 63    & 31    & 31    & 6     & 1.7487 & 1,044 & 597   & 6     & 10    & 99.5  & 6.8 \\
A-n63-k9b & 63    & 46    & 16    & 8     & 1.7898 & 1,405 & 785   & 8     & 8     & 98.1  & 6.8 \\
A-n64-k9a & 64    & 31    & 32    & 6     & 1.5874 & 935   & 589   & 6     & 9     & 98.2  & 6.7 \\
A-n64-k9b & 64    & 47    & 16    & 7     & 1.6440 & 1,136 & 691   & 7     & 5     & 98.7  & 7.4 \\
A-n65-k9a & 65    & 32    & 32    & 6     & 1.3232 & 786   & 594   & 6     & 10    & 99.0  & 7.0 \\
A-n65-k9b & 65    & 48    & 16    & 8     & 1.3197 & 1,036 & 785   & 8     & 9     & 98.1  & 7.1 \\
A-n69-k9a & 69    & 34    & 34    & 5     & 1.4061 & 696   & 495   & 5     & 10    & 99.0  & 8.8 \\
A-n69-k9b & 69    & 51    & 17    & 7     & 1.3968 & 968   & 693   & 7     & 7     & 99.0  & 8.3 \\
A-n80-k10a & 80    & 39    & 40    & 7     & 1.8371 & 1,286 & 700   & 7     & 14    & 100.0 & 7.6 \\
A-n80-k10b & 80    & 59    & 20    & 9     & 1.8401 & 1,623 & 882   & 9     & 11    & 98.0  & 7.8 \bigstrut[b]\\
\hline
\end{tabular}%

}}
\end{center}
\end{table}

\begin{table}
\caption{Details of the best solutions found for class \classPA} \label{tab:sol-PA}
\begin{center}
{\scriptsize{
\setlength{\tabcolsep}{2.7pt}
\renewcommand{\arraystretch}{1.0}
\begin{tabular}{rrrrr|rrrrrrr}
\hline
\tabname & $\setv$ & $\nf$ & $\nc$ & $\nveh$ & \tabzbest & Profit  & Time  & \tabsolm & \tabnfixed & \tabut & \tabavgc \bigstrut\\
\hline
A-n32-k5a & 32    & 15    & 16    & 3     & 0.3835 & 306   & 798   & 3     & 5     & 97.8  & 6.7 \bigstrut[t]\\
A-n32-k5b & 32    & 23    & 8     & 4     & 0.3809 & 406   & 1,066 & 4     & 6     & 98.0  & 7.3 \\
A-n33-k5a & 33    & 16    & 16    & 3     & 0.6014 & 433   & 720   & 3     & 13    & 98.0  & 9.7 \\
A-n33-k5b & 33    & 24    & 8     & 4     & 0.5989 & 433   & 723   & 3     & 5     & 98.4  & 9.7 \\
A-n33-k6a & 33    & 16    & 16    & 3     & 0.7856 & 458   & 583   & 3     & 10    & 82.3  & 8.7 \\
A-n33-k6b & 33    & 24    & 8     & 4     & 0.6521 & 506   & 776   & 4     & 5     & 81.5  & 7.3 \\
A-n34-k5a & 34    & 16    & 17    & 3     & 0.5759 & 368   & 639   & 3     & 9     & 91.4  & 8.3 \\
A-n34-k5b & 34    & 24    & 9     & 4     & 0.5617 & 437   & 778   & 4     & 5     & 83.1  & 7.3 \\
A-n36-k5a & 36    & 17    & 18    & 4     & 0.3834 & 403   & 1,051 & 4     & 13    & 99.9  & 7.5 \\
A-n36-k5b & 36    & 26    & 9     & 4     & 0.4044 & 421   & 1,041 & 4     & 6     & 99.0  & 8.0 \\
A-n37-k5a & 37    & 18    & 18    & 3     & 0.4067 & 305   & 750   & 3     & 5     & 98.4  & 7.7 \\
A-n37-k5b & 37    & 27    & 9     & 5     & 0.4325 & 397   & 918   & 4     & 7     & 90.4  & 8.5 \\
A-n37-k6a & 37    & 18    & 18    & 4     & 0.5044 & 516   & 1,023 & 4     & 11    & 98.7  & 7.3 \\
A-n37-k6b & 37    & 27    & 9     & 6     & 0.4371 & 559   & 1,279 & 5     & 7     & 98.8  & 6.8 \\
A-n38-k5a & 38    & 18    & 19    & 4     & 0.5111 & 437   & 855   & 4     & 12    & 86.2  & 7.5 \\
A-n38-k5b & 38    & 27    & 10    & 4     & 0.4922 & 476   & 967   & 4     & 8     & 97.5  & 8.8 \\
A-n39-k5a & 39    & 19    & 19    & 4     & 0.5070 & 397   & 783   & 4     & 8     & 78.6  & 6.8 \\
A-n39-k5b & 39    & 28    & 10    & 5     & 0.4751 & 449   & 945   & 4     & 4     & 94.9  & 8.0 \\
A-n39-k6a & 39    & 19    & 19    & 4     & 0.5310 & 403   & 759   & 3     & 4     & 95.8  & 7.7 \\
A-n39-k6b & 39    & 28    & 10    & 5     & 0.4747 & 479   & 1,009 & 4     & 4     & 95.5  & 8.0 \\
A-n44-k6a & 44    & 21    & 22    & 4     & 0.4084 & 437   & 1,070 & 4     & 10    & 98.7  & 7.8 \\
A-n44-k6b & 44    & 32    & 11    & 5     & 0.4470 & 540   & 1,208 & 5     & 7     & 89.2  & 7.8 \\
A-n45-k6a & 45    & 22    & 22    & 5     & 0.5667 & 514   & 907   & 4     & 14    & 80.7  & 9.0 \\
A-n45-k6b & 45    & 33    & 11    & 5     & 0.5261 & 585   & 1,112 & 4     & 9     & 98.6  & 10.5 \\
A-n45-k7a & 45    & 22    & 22    & 5     & 0.4049 & 462   & 1,141 & 5     & 10    & 98.8  & 6.4 \\
A-n45-k7b & 45    & 33    & 11    & 7     & 0.3666 & 569   & 1,552 & 7     & 5     & 96.0  & 5.4 \\
A-n46-k7a & 46    & 22    & 23    & 4     & 0.5554 & 531   & 956   & 4     & 13    & 98.4  & 8.8 \\
A-n46-k7b & 46    & 33    & 12    & 5     & 0.5672 & 574   & 1,012 & 5     & 7     & 83.3  & 8.0 \\
A-n48-k7a & 48    & 23    & 24    & 5     & 0.4369 & 530   & 1,213 & 5     & 12    & 95.5  & 7.0 \\
A-n48-k7b & 48    & 35    & 12    & 6     & 0.3955 & 598   & 1,512 & 6     & 9     & 99.2  & 7.3 \\
A-n53-k7a & 53    & 26    & 26    & 4     & 0.5405 & 554   & 1,025 & 4     & 13    & 98.9  & 9.8 \\
A-n53-k7b & 53    & 39    & 13    & 7     & 0.4699 & 594   & 1,264 & 5     & 5     & 95.0  & 8.8 \\
A-n54-k7a & 54    & 26    & 27    & 5     & 0.4814 & 594   & 1,234 & 5     & 15    & 98.7  & 8.2 \\
A-n54-k7b & 54    & 39    & 14    & 7     & 0.3949 & 605   & 1,532 & 7     & 7     & 87.2  & 6.6 \\
A-n55-k9a & 55    & 27    & 27    & 5     & 0.7779 & 746   & 959   & 4     & 16    & 94.0  & 10.8 \\
A-n55-k9b & 55    & 40    & 14    & 5     & 0.8000 & 760   & 950   & 4     & 8     & 93.1  & 12.0 \\
A-n60-k9a & 60    & 29    & 30    & 5     & 0.5878 & 656   & 1,116 & 5     & 13    & 85.8  & 8.4 \\
A-n60-k9b & 60    & 44    & 15    & 6     & 0.5281 & 760   & 1,439 & 6     & 8     & 92.2  & 8.7 \\
A-n61-k9a & 61    & 30    & 30    & 5     & 0.7976 & 717   & 899   & 4     & 16    & 97.3  & 11.5 \\
A-n61-k9b & 61    & 45    & 15    & 6     & 0.7754 & 825   & 1,064 & 5     & 9     & 92.1  & 10.8 \\
A-n62-k8a & 62    & 30    & 31    & 5     & 0.4335 & 619   & 1,428 & 5     & 16    & 96.5  & 9.2 \\
A-n62-k8b & 62    & 45    & 16    & 7     & 0.4085 & 719   & 1,760 & 6     & 11    & 92.5  & 9.3 \\
A-n63-k10a & 63    & 31    & 31    & 5     & 0.6993 & 800   & 1,144 & 5     & 16    & 89.7  & 9.4 \\
A-n63-k10b & 63    & 46    & 16    & 7     & 0.6105 & 887   & 1,453 & 6     & 9     & 94.6  & 9.2 \\
A-n63-k9a & 63    & 31    & 31    & 7     & 0.3203 & 583   & 1,820 & 7     & 7     & 99.2  & 5.4 \\
A-n63-k9b & 63    & 46    & 16    & 9     & 0.3597 & 838   & 2,330 & 9     & 10    & 98.8  & 6.2 \\
A-n64-k9a & 64    & 31    & 32    & 7     & 0.4333 & 770   & 1,777 & 7     & 19    & 98.8  & 7.1 \\
A-n64-k9b & 64    & 47    & 16    & 6     & 0.4548 & 821   & 1,805 & 6     & 12    & 97.4  & 9.8 \\
A-n65-k9a & 65    & 32    & 32    & 7     & 0.6879 & 787   & 1,144 & 5     & 20    & 97.4  & 10.4 \\
A-n65-k9b & 65    & 48    & 16    & 8     & 0.6189 & 851   & 1,375 & 6     & 13    & 97.1  & 10.2 \\
A-n69-k9a & 69    & 34    & 34    & 4     & 0.6641 & 678   & 1,021 & 4     & 16    & 99.7  & 12.5 \\
A-n69-k9b & 69    & 51    & 17    & 6     & 0.5939 & 756   & 1,273 & 5     & 8     & 99.1  & 11.8 \\
A-n80-k10a & 80    & 39    & 40    & 6     & 0.3911 & 761   & 1,946 & 6     & 16    & 99.8  & 9.2 \\
A-n80-k10b & 80    & 59    & 20    & 7     & 0.3894 & 903   & 2,319 & 7     & 11    & 98.9  & 10.0 \bigstrut[b]\\
\hline
\end{tabular}%

}}
\end{center}
\end{table}


\end{document}